\documentclass[11pt]{article}
\usepackage{amsfonts,amssymb,amsmath,amsthm}
\usepackage{url}
\usepackage{enumerate}

\usepackage{dsfont,mathrsfs,stmaryrd}
\input xy
\xyoption{all}

\newtheorem{thm}{Theorem}[subsection]
\newtheorem{cor}[thm]{Corollary}
\newtheorem{lem}[thm]{Lemma}
\newtheorem{prop}[thm]{Proposition}
\newtheorem*{thm2}{Theorem}

\theoremstyle{plain}
\newtheorem{thmm}{Theorem}[section]

\newtheorem{lem2}[thmm]{Lemma}

\theoremstyle{definition}
\newtheorem{defn}[thm]{Definition}

\newtheorem*{remark}{\textbf{Remark}}

\newcommand{\bd}{\begin{defn}}
\newcommand{\ed}{\end{defn}}
\newcommand{\bl}{\begin{lem}}
\newcommand{\el}{\end{lem}}
\newcommand{\bp}{\begin{prop}}
\newcommand{\ep}{\end{prop}}
\newcommand{\bt}{\begin{thm}}
\newcommand{\et}{\end{thm}}
\newcommand{\bc}{\begin{cor}}
\newcommand{\ec}{\end{cor}}
\newcommand{\br}{\begin{remark}}
\newcommand{\er}{\end{remark}}
\newcommand{\bdi}{\begin{diagram}}
\newcommand{\edi}{\end{diagram}}
\newcommand{\beq}{\begin{eqn}}
\newcommand{\eeq}{\end{eqn}}
\newcommand{\ba}{\begin{array}}
\newcommand{\ea}{\end{array}}
\newcommand{\bpf}{\begin{proof}}
\newcommand{\epf}{\end{proof}}

\newcommand{\pf}{\noindent {\it Proof}}

\newcommand{\Z}{\mathds{Z}}
\newcommand{\Q}{\mathds{Q}}
\newcommand{\Zp}{\mathds{Z}_{p}}
\newcommand{\Qp}{\mathds{Q}_{p}}
\newcommand{\al}{\alpha}
\newcommand{\be}{\beta}
\newcommand{\de}{\delta}

\newcommand{\ga}{\gamma}
\newcommand{\Ga}{\Gamma}

\newcommand{\e}{\varepsilon}

\newcommand{\s}{\sigma}

\newcommand{\La}{\Lambda}
\newcommand{\la}{\lambda}

\newcommand{\F}{\mathscr{F}}
\newcommand{\U}{\mathscr{U}}
\newcommand{\C}{\mathcal{C}}
\newcommand{\Di}{\mathcal{D}}

\newcommand{\ai}{\mathfrak{a}}

\newcommand{\Gal}{\mathrm{Gal}}
\newcommand{\Hom}{\mathrm{Hom}}

\newcommand{\Cone}{\mathrm{Cone}}
\newcommand{\ot}{\otimes}
\newcommand{\otw}{\hat{\ot}}
\newcommand{\bo}{\bigoplus}

\newcommand{\ilim}{\displaystyle \mathop{\varinjlim}\limits}
\newcommand{\plim}{\displaystyle \mathop{\varprojlim}\limits}

\newcommand{\cts}{\mathrm{cts}}

\newcommand{\cd}{\mathrm{cd}}
\newcommand{\DD}{\mathbf{D}}
\newcommand{\R}{\mathbf{R}}

\newcommand{\ra}{\rightarrow}
\newcommand{\lra}{\longrightarrow}
\newcommand{\tha}{\twoheadrightarrow}

\newcommand{\hra}{\hookrightarrow}

\newcommand{\sbs}{\subseteq}

\begin{document}

\title{Poitou-Tate duality over extensions of global fields}

\author{Meng Fai Lim}
\date{}
\maketitle

\begin{abstract}
\noindent In this paper, we are interested in the Poitou-Tate
duality in Galois cohomology. We will formulate and prove a theorem
for a nice class of modules (with a continuous Galois action) over a
pro-$p$ ring. The theorem will comprise of the Tate local duality,
Poitou-Tate duality and the Poitou-Tate's exact sequence.
    \end{abstract}

\section{Introduction}

The classical Poitou-Tate duality is a duality principle for a
local-global statement, namely it relates the kernels of the
localization maps. Using compactly supported cohomology groups, one
can give a cleaner formulation of the statement which we now do. Let
$F$ be a global field with characteristic not equal to $p$, and let
$S$ be a finite set of primes of $F$ containing all primes above $p$
and all archimedean primes of $F$. We let $G_{F,S}$ denote the
Galois group $\Gal(F_S/F)$ of the maximal unramified outside $S$
extension $F_S$ of $F$ inside a fixed separable closure of $F$. In
its usual formulation, Poitou-Tate duality relates the kernels of
the localization maps on the $G_{F,S}$-cohomology of a module and
the Tate twist of its Pontryagin dual. For simplicity, we assume in
this introduction that $p$ is odd if $F$ has any real places. The
general result without this assumption can be found in Theorem
\ref{Poitou-Tate adic}.

The $n$th compactly supported $G_{F,S}$-cohomology group
$H^n_{c,\cts}(G_{F,S},M)$ with coefficients in a topological
$G_{F,S}$-module $M$ is defined as the $n$th cohomology group of the
complex
\[ \Cone\left( C^{\cdot}_{\cts}(G_{F,S},M) \stackrel{\mathrm{res}_S}{\lra}
 \bigoplus_{v \in S_f} C^{\cdot}_{\cts}(G_{F_v},M)\right)[-1],
\]
where $G_{F_v}$ is the absolute Galois group of the completion of
$F$ at $v$, and $\mathrm{res}_S$ is the sum of restriction maps on
the continuous cochain complexes. It therefore fits in a long exact
sequence
$$
    \cdots \to H^n_{c,\cts}(G_{F,S},M) \to H^n_{\cts}(G_{F,S},M) \to \bigoplus_{v \in S} H^n_{\cts}
    (G_{F_v},M) \to H^{n+1}_{c,\cts}(G_{F,S},M) \to \cdots.
$$

We now let $R$ denote a commutative complete Noetherian local ring
with finite residue field of characteristic $p$.  Then we have the
following formulation of Poitou-Tate duality due to
Nekov$\mathrm{\acute{a}}\mathrm{\check{r}}$ \cite[Prop.\
5.4.3(i)]{Ne}.

\begin{thm2}
Let $T$ be a finitely generated $R$-module with a continuous
$(R$-linear$)$ $G_{F,S}$-action. Then there are isomorphisms
\[ \ba{c} H^n_{\cts}(G_{F,S},T) \stackrel{\sim}{\lra} H^{3-n}_{c, \cts}(G_{F,S}, T^{\vee}(1))^{\vee} \\
     H^n_{c,\cts}(G_{F,S},T) \stackrel{\sim}{\lra} H^{3-n}_{\cts}(G_{F,S}, T^{\vee}(1))^{\vee} \\
\ea\] of $R$-modules for all $n$, where $T^{\vee} = \Hom_{\cts}(T,
\Qp/\Zp)$.
\end{thm2}

We now recall some notation from the language of derived categories.
We denote by $\DD(Mod_R)$ the derived category of $R$-modules which
is obtained from from the category $\mathrm{Ch}(Mod_R)$ of chain
complexes of $R$-modules by inverting the quasi-isomorphisms, i.e.,
the maps of complexes that induce isomorphisms on cohomology. We
have the derived functors $\R\Hom_R(-,-)$, $\R\Ga_{\cts}(G_{F,S},-)$
and $\R\Ga_{c,\cts}(G_{F,S},-)$ that are obtained from
$\Hom_R(-,-)$, $C_{\cts}(G_{F,S},-)$ and $C_{c, \cts}(G_{F,S},-)$.
Then Poitou-Tate duality can be reformulated as the following
isomorphisms \[ \ba{c} \R\Ga_{\cts}(G_{F,S},T) \stackrel{\sim}{\lra}
\R\Hom_{\Zp}\Big(\R\Ga_{c, \cts}(G_{F,S}, T^{\vee}(1)),
\Qp/\Zp\Big)[-3] \\
 \R\Ga_{c,\cts}(G_{F,S},T) \stackrel{\sim}{\lra}
\R\Hom_{\Zp}\Big(\R\Ga_{\cts}(G_{F,S}, T^{\vee}(1)),
\Qp/\Zp\Big)[-3] \ea
\]
in $\DD(Mod_R)$.

Now suppose that $F_{\infty}$ is a $p$-adic Lie extension of $F$
contained in $F_S$. We denote by $\Ga$ the Galois group of the
extension $F_{\infty}/F$, and we let $\Lambda =
R\llbracket\Gamma\rrbracket$ denote the resulting Iwasawa algebra
over $R$. Let $T$ be a finitely generated $R$-module with a
continuous ($R$-linear) $G_{F,S}$-action, and let $A$ be a
cofinitely generated $R$-module with a continuous ($R$-linear)
$G_{F,S}$-action. The $\Lambda$-modules of interest are the
following direct and inverse limits of cohomology groups (and their
counterparts with compact support)
\[ \ba{rcl} \ilim_{F_{\al}} H^n_{\cts}(\mathrm{Gal}(F_S/F_{\al}),A) &\text{and}& \plim_{F_{\al}}
 H^n_{\cts}(\mathrm{Gal}(F_S/F_{\al}), T),\ea \] where the limits are taken over all
finite Galois extensions $F_{\al}$ of $F$ which are contained in
$F_{\infty}$. By an application of Shapiro's lemma, one can show
that they are respectively isomorphic to
\[ \ba{rcl} H^n_{\cts}(G_{F,S}, F_{\Ga}(A)) &\text{and}&
H^n_{\cts}(G_{F,S}, \F_{\Ga}(T)), \ea\] where the $\La$-modules
$F_{\Ga}(A)$ and $\F_{\Ga}(T)$ are defined by
\[ \ba{rcl} \ilim_{F_{\al}}\Hom_R(R[\mathrm{Gal}(F_{\al}/F)], A) &\text{and}&
\plim_{F_{\al}} R[\mathrm{Gal}(F_{\al}/F)]\ot_R T \ea\]
respectively. Therefore, we can reduce the question of finding
dualities on the Iwasawa modules of interest to that of obtaining
dualities over $G_{F,S}$, but with $R$ replaced by $\La$.

In his monograph \cite{Ne},
Nekov$\mathrm{\acute{a}}\mathrm{\check{r}}$ considers the above
situation over a commutative $p$-adic Lie extension (e.g., a
$\Zp^r$-extension) and develops an extension of Poitou-Tate global
duality for the above cohomology groups. In view of the vast
activity in the study of noncommutative generalizations of the main
conjecture of Iwasawa theory \cite{CFKSV, FK, Ka, RW}, one would
like to extend the above theory to the noncommutative setting.

In fact, in this paper, we study generalizations of the above
duality of Poitou-Tate over a general pro-$p$ ring $\La$ (not
necessarily commutative). Together with the module theory, we
carefully develop the theory of continuous group cohomology in our
setting. From there, we are able to state and prove our duality
theorem (cf. Theorem \ref{Poitou-Tate adic}).

\begin{thm2} Let $M$ be a bounded
complex of objects that are profinite $\La$-modules with a
continuous $(\La$-linear$)$ $G_{F,S}$-action. Then we have the
following isomorphism
\[ \entrymodifiers={!! <0pt, .8ex>+} \xymatrix @C=0.7in{
  \displaystyle\bigoplus_{v\in S}\R\Ga(G_{v}, M)[-1] \ar[d]_{} \ar[r]^{} &
  \displaystyle\bigoplus_{v\in S}\R\Hom_{\Zp}\big(\R\Ga(G_{v}, M^{\vee}(1)), \Qp/\Zp\big)[-3] \ar[d]  \\
  \R\Ga_c(G_{F,S}, M) \ar[d]_{} \ar[r]^{} &
  \R\Hom_{\Zp}\big(\R\Ga(G_{F,S}, M^{\vee}(1)), \Qp/\Zp\big)[-3] \ar[d]  \\
  \R\Ga(G_{F,S}, M) \ar[r]^{} & \R\Hom_{\Zp}\big(\R\Ga_c(G_{F,S}, M^{\vee}(1)), \Qp/\Zp\big)[-3]   }
  \] of exact triangles in $\DD(Mod_{\La})$. \end{thm2}

We now give a brief description of the contents of each section of
the paper. In Section \ref{Preliminaries}, we introduce notations
and results from homological algebra required for the paper. Section
\ref{Profinite rings} is about the discussion of profinite rings and
their topological modules. We also introduce continuous cohomology
groups with coefficients in compact modules and discrete modules. In
Section \ref{duality}, we will formulate and prove our duality
theorems. In Section \ref{Iwasawa modules}, we will apply the
duality theorems proved in Section \ref{duality} to extensions of
global fields.

\proof[Acknowledgements]
    The material presented in this article
    forms a generalized version of first part of the author's Ph.D.\ thesis
    \cite{Lim-thesis}. The remaining part of the thesis,
    for which this paper provides some preliminaries, can be found in \cite{LS}.  The author would like to thank his
    Ph.D.\ supervisor Romyar Sharifi for his advice and
    encouragement, for without which, this paper would not have been
    possible. The author would also like to thank Manfred Kolster for his
    encouragement. Many thanks also go to the referee for a number of
    comments and suggestions.

\section{Preliminaries} \label{Preliminaries}

We begin by reviewing certain objects and notation that will be used
in this write-up. Most of the material presented in this section can
be found in \cite{RD, Ne, Wei}. Throughout the paper, every ring is
associative and has a unit.

Fix an abelian category $\mathcal{A}$ and denote the category of
(cochain) complexes of objects in $\mathcal{A}$ by
Ch($\mathcal{A}$). We also denote the category of bounded below
complexes, bounded above complexes and bounded complexes by
Ch$^+$($\mathcal{A}$), Ch$^-(\mathcal{A}$) and Ch$^b(\mathcal{A}$)
respectively. For each $n\in\Z$, the translation by $n$ of a complex
$X$ is given by \[ X[n]^{i} = X^{n+i}, ~d^{i}_{X[n]} =
(-1)^{n}d^{n+i}_{X}.\] If $f:X\lra Y$ is a morphism of complexes,
then $f[n]:X[n]\lra Y[n]$ is given by $f[n]^{i} = f^{n+i}$.

If $X$ is a complex, we have the following truncations of $X$:
\[ \ba{l}
\s_{\leq i}X = [ \cdots\lra X^{i-2}\lra X^{i-1}\lra X^i \lra 0 \lra
0 \lra \cdots ] \\
 \tau_{\leq i}X = [ \cdots\lra X^{i-2}\lra X^{i-1}\lra \ker(d_{X}^i) \lra 0 \lra
0 \lra \cdots ] \\
 \s_{\geq i}X = [ \cdots\lra 0\lra 0\lra X^i \lra X^{i+1} \lra
X^{i+2} \lra \cdots ]\\
 \tau_{\geq i}X = [ \cdots\lra 0\lra 0\lra \mathrm{coker}(d_{X}^{i-1}) \lra X^{i+1} \lra
X^{i+2} \lra \cdots ].\ea \]

The cone of a morphism $f:X\lra Y$ is defined by $\mathrm{Cone}(f) =
Y\oplus X[1]$ with differential
\[ d^{i}_{\mathrm{Cone}(f)} = \left(\begin{array}{cc}
                                             d^{i}_{Y} & f^{i+1} \\
                                             0 & -d^{i+1}_{X} \\
                                           \end{array}\right) :
Y^{i}\oplus X^{i+1} \lra Y^{i+1}\oplus X^{i+2}. \] There is an exact
sequence of complexes
\[ 0\lra Y\stackrel{j}{\lra} \Cone(f)\stackrel{p}{\lra} X[1]\lra 0,\]
where $j$ and $p$ are the canonical inclusion and projection
respectively. The corresponding boundary map
\[ \delta : H^{i}(X[1]) = H^{i+1}(X) \lra H^{i+1}(Y) \]
is induced by $f^{i+1}$.

If $X$ is a complex and $x\in X^i$, we write $\bar{x} = i$ for the
degree.

Let $\La, S$ and $T$ be rings. Let $M$ (resp., $N$) be a left
$\La$-$S$-bimodule (resp., a $\La$-$T$-bimodule). Then
$\Hom_{\La}(M,N)$ is taken to be the $S$-$T$-bimodule of all left
$\La$-module homomorphisms from $M$ to $N$, where the left
$S$-action is given by $(s\cdot f)(m) = f(ms)$ and the right
$T$-action is given by $(f\cdot t)(m) = f(m)t$ for $f\in
\Hom_{\La}(M,N), m\in M, s\in S$ and $t\in T$. If $M^{\bullet}$ is a
complex of $\La$-$S$-bimodules and $N^{\bullet}$ a complex of
$\La$-$T$-bimodules, we define a complex
$\Hom_{\La}^{\bullet}(M^{\bullet},N^{\bullet})$ of $S$-$T$-bimodules
by
\[ \Hom_{\La}^{n}(M^{\bullet},N^{\bullet}) = \prod_{i\in\Z}\Hom_{\La}(M^{i}, N^{i+n}) \]
with differentials defined as follows: for $f\in
\Hom_{\La}(M^i,N^{i+n})$, we have
\[ df = d^{i+n}_N\circ f + (-1)^{n}f\circ d^{i-1}_M .\]

In the case when $S=T$, we have a similar definition for the
complexes $\Hom^{\bullet}_{\La-S}(M^{\bullet}, N^{\bullet})$ of
abelian groups, where $\Hom_{\La-S}(M,N)$ is the group of all of
$\La$-$S$-bimodule homomorphisms from $M$ to $N$. It follows
immediately from the definition that for an element $f \in
\Hom^0_{\La-S}(M^{\bullet}, N^{\bullet})$, we have $f\in
\Hom_{\mathrm{Ch}(\La-S)}(M^{\bullet}, N^{\bullet})$ if and only if
$df = 0$. Here $\mathrm{Ch}(\La-S)$ denotes the category of
complexes of $\La$-$S$-bimodules.

Suppose that $M^{\bullet}$ is a complex of $\La$-$S$-bimodules and
$L^{\bullet}$ a complex of $S$-$T$-bimodules. We define the complex
$M^{\bullet}\ot_{S}L^{\bullet}$ of $\La$-$T$-bimodules by
\[ (M^{\bullet}\ot_{S}L^{\bullet})^{n} = \bigoplus_{i\in\Z}M^{i}\ot_{S}L^{n-i} \]
with differentials
\[ d(m\ot l) = dm\ot l + (-1)^{\bar{m}} m\ot dl .\]

We end the section by collecting some technical results which will
be used in the paper.

\begin{lem2} \label{translation lemma}The following formulas define isomorphisms of complexes$:$
\[
\ba{c} \Hom_{\La}^{\bullet}(M^{\bullet},N^{\bullet})[n] \cong
\Hom_{\La}^{\bullet}(M^{\bullet},N^{\bullet}[n]) \\
 f\mapsto f \\
(M^{\bullet}[n])\ot_S L^{\bullet}\cong (M^{\bullet}\ot_S
L^{\bullet})[n]\\
m\ot l \mapsto m\ot l\\
M^{\bullet}\ot_S (L^{\bullet}[n])\cong (M^{\bullet}\ot_S
L^{\bullet})[n]\\
m\ot l \mapsto (-1)^{n\bar{m}}m\ot l.\\ \ea \]
 \end{lem2}

\bpf This follows from a straightforward verification of the
definition of translation and the sign conventions. \epf

\begin{lem2} \label{adj sign} The adjunction morphisms define morphisms
\[ \ba{c} \Hom^{\bullet}_{\La-T}(M^{\bullet}\ot_S L^{\bullet}, N^{\bullet})
\lra \Hom^{\bullet}_{\La-S}\big(M^{\bullet},
\Hom^{\bullet}_{T^o}(L^{\bullet}, N^{\bullet})\big) \\
 f\mapsto \big( m \mapsto(l\mapsto f(m\ot l))\big) \\
 \Hom^{\bullet}_{\La-T}(M^{\bullet}\ot_S L^{\bullet}, N^{\bullet})
\lra \Hom^{\bullet}_{S-T}\big(L^{\bullet},
\Hom^{\bullet}_{\La}(M^{\bullet}, N^{\bullet})\big) \\
 f\mapsto \big( l \mapsto(m\mapsto (-1)^{\bar{m}\bar{l}} f(m\ot
 l))\big)\ea \] of complexes and morphisms
\[ \ba{c} \Hom_{\mathrm{Ch}(\La-T)}(M^{\bullet}\ot_S L^{\bullet}, N^{\bullet})
\lra \Hom_{\mathrm{Ch}(\La-S)}\big(M^{\bullet},
\Hom^{\bullet}_{T^o}(L^{\bullet}, N^{\bullet})\big) \\
\Hom_{\mathrm{Ch}(\La-T)}(M^{\bullet}\ot_S L^{\bullet}, N^{\bullet})
\lra \Hom_{\mathrm{Ch}(S-T)}\big(L^{\bullet},
\Hom^{\bullet}_{\La}(M^{\bullet}, N^{\bullet})\big) \ea
\] of abelian groups. All of these maps are monomorphisms; they are isomorphisms if
$M^{\bullet}$ and $L^{\bullet}$ are bounded above and $N^{\bullet}$
is bounded below.  \end{lem2}

\begin{lem2} \label{cone tech} Given the following data:

 $(1)$ Complexes $A_1$, $B_1$ of $\La$-$S$-bimodules, complexes $A_2$, $B_2$ of
 $S$-$T$-bimodules, and complexes $A_3, B_3$ of $\La$-$T$-bimodules.

 $(2)$ Morphisms of complexes $f_j :A_j \lra B_j$ preserving the
 respective bimodule structures.

 $(3)$ Morphisms of complexes of $\La$-$T$-bimodules
 \[ \ba{c} \cup_A : A_1\ot_S A_2 \lra A_3 \\
            \cup_B : B_1\ot_S B_2 \lra B_3 \ea \] such that $f_3\circ
            \cup_A = \cup_B \circ(f_1\ot f_2)$.
For $j=1,2,3$, define $E_j$ to be the complex
\[ \Cone\big( A_j\stackrel{f_j}{\lra} B_j\big)[-1]. \]
Then we have morphisms of complexes
\[\cup_0, \cup_1 : E_1\ot_S E_2\lra E_3\] given by the formulas
\[ \ba{c} (a_1,b_1)\cup_0(a_2, b_2) = \big( a_1\cup_A a_2,
(-1)^{\bar{a}_1}f_1(a_1)\cup_B b_2\big) \\
     (a_1,b_1)\cup_1(a_2, b_2) = \big( a_1\cup_A a_2,
b_1\cup_B f_2(a_2)\big) ,           \ea \] and the formula

\[ s\big((a_1,b_1)\ot(a_2,b_2)\big) = \big(0,
(-1)^{\bar{a}_1}b_1 \cup_B b_2\big)\] defines a homotopy $s :
\cup_1\rightsquigarrow \cup_0$. \end{lem2}

\bpf This is a special case of \cite[Prop.\ 1.3.2]{Ne}. \epf

\section{Profinite rings} \label{Profinite rings}

Completed group algebras of profinite groups arise naturally in the
study of Iwasawa theory, and such rings are profinite rings. In this
section, we shall study the properties of profinite rings and their
(topological) modules. We will also develop a cohomological theory
over such rings.

Throughout the section, $\La$ will always denote a profinite ring,
and $\mathcal{I}$ is a directed fundamental system of open
neighborhoods of zero consisting of two-sided ideals of $\La$. We
use $\La^{\circ}$ to denote the opposite ring to $\La$.

\subsection{Topological $\La$-modules}

In this subsection, we will study the topological modules over a
profinite ring $\La$. These are Hausdorff topological abelian groups
with a continuous $\La$-action. In particular, we are interested in
the following two classes of topological $\La$-modules.

\bd We say that a topological $\La$-module $M$ is a compact (resp.,
discrete) $\La$-module if its underlying topology is compact (resp.,
discrete). The category of compact $\La$-modules (resp., discrete
$\La$-modules) is denoted by $\C_{\La}$ (resp., $\Di_{\La}$).\ed

The following proposition records some of the properties of the
above two categories, whose proofs can be found in \cite[Chap.\
5]{RZ}.

\bp \label{Pontryagin duality} $(i)$ Every compact $\La$-module is a
projective limit of finite modules and has a fundamental system of
neighborhoods of zero consisting of open submodules. In particular,
it is an abelian profinite group.

$(ii)$ Every discrete $\La$-module is the direct limit of finite
$\La$-modules. In particular, it is an abelian torsion group.

$(iii)$ Pontryagin duality induces a duality between the category
$\C_{\La}$ of compact $\La$-modules and the category
$\Di_{\La^{\circ}}$ of discrete $\La^{\circ}$-modules.

$(iv)$ The category $\C_{\La}$ is abelian and has enough projectives
and exact inverse limits. The category $\Di_{\La}$ is abelian and
has enough injectives and exact direct limits. \ep

\medskip We give another description of discrete $\La$-modules. If
$M$ is a $\La$-module and $\ai$ is a two-sided ideal of $\La$, we
define
\[ M[\ai] = \{ x\in M\,  |~ \ai\sbs \mathrm{Ann}(x)\}.\] With
this, we have the following lemma.

\bl \label{discrete La} Let $M$ be an abstract $\La$-module. Then
$M$ is a discrete $\La$-module $($i.e.,\ the $\La$-action is
continuous with respect to the discrete topology on $M$$)$ if and
only if
\[ M = \displaystyle \bigcup_{\ai\in\mathcal{I}}M[\ai]. \] \el

\bpf Suppose that $M$ is a discrete $\La$-module. Let $x\in M$. Then
by the continuity of the $\La$-action, there exists
$\ai\in\mathcal{I}$ such that $\ai\cdot x = 0$. This implies that
$x\in M[\ai]$.

Conversely, suppose that \[ M = \displaystyle
\bigcup_{\ai\in\mathcal{I}}M[\ai]. \] We shall show that the action
\[ \theta :\La\times M \lra M \] is continuous, where $M$ is given
the discrete topology. In other words, for each $x\in M$, we need to
show that $\theta^{-1}(x)$ is open in $\La\times M$. Let $(\la,
y)\in \theta^{-1}(x)$. Then $y\in M[\ai]$ for some
$\ai\in\mathcal{I}$. Therefore, we have $(\la, y) \in
(\la+\ai)\times\{y\}$, and the latter set is an open set contained
in $\theta^{-1}(x)$. \epf

\medskip When working with topological $\La$-modules, one will have
to consider continuous homomorphisms between the modules. In
general, an abstract homomorphism of modules may not be continuous.
In the next lemma, we record a few situations where every abstract
homomorphism is continuous. We say that a topological $\La$-module
$M$ is endowed with the $\mathcal{I}$-adic topology if the
collection $\{\ai M\}_{\ai\in\mathcal{I}}$ forms a fundamental
system of neighborhoods of zero.

 \bl \label{abstract is cts} Let $M$ and $N$
be two topological $\La$-modules. Suppose one of the following cases
holds.

$(1)$ Both $M$ and $N$ have the $\mathcal{I}$-adic topology.

$(2)$ Both $M$ and $N$ have the discrete topology.

$(3)$ $M$ is a finitely generated $\La$-module endowed with the
 $\mathcal{I}$-adic topology, and $N$ is a compact $\La$-module.

$(4)$ $M$ is a finitely generated $\La$-module endowed with the
 $\mathcal{I}$-adic topology, and $N$ is a discrete $\La$-module.

Then every abstract $\La$-homomorphism is continuous. In other
words, we have
\[ \Hom_{\La,\,\cts}(M,N) = \Hom_{\La}(M,N). \] \el

\bpf (1) and (2) are straightforward.

(3) Suppose $M$ is generated by $e_1,..., e_r$.  Let $f:M\lra N$ be
an abstract $\La$-homomorphism, and for each $i$, set $x_i =
f(e_i)$. Let $V$ be an open $\La$-submodule of $N$. By continuity of
the $\La$-action on $N$, for each $i$, there exists $\ai_i\in
\mathcal{I}$ such that $\ai_i\cdot x_i\sbs V$. Since $\mathcal{I}$
is directed, we can find $\ai\in\mathcal{I}$ such that $\ai\sbs
\ai_i$ for all $i$. It follows that $f(\ai M) \sbs V$, establishing
the continuity of $f$.

(4) We retain the notation in (3). By Lemma \ref{discrete La}, for
each $i$, there exists $\ai_i\in\mathcal{I}$ such that $x_i\in
N[\ai_i]$. Since $\mathcal{I}$ is directed, we can find
$\ai\in\mathcal{I}$ such that $\ai\sbs \ai_i$ for all $i$, and
$f(\ai M) =0$. \epf

\bc \label{fg submodule is closed} Let $M$ be a compact
$\La$-module. Then every finitely generated abstract $\La$-submodule
of $M$ is a closed subset of $M$. In particular, every finitely
generated left $($or right$)$ ideal of $\La$ is closed in $\La$. \ec

\bpf  Let $N$ be a $\La$-submodule of $M$ generated by
$x_{1},...,x_{r}$. By Lemma \ref{abstract is cts}(3), the following
$\La$-homomorphism
\[ \ba{rl}
   \phi:\bigoplus_{i=1}^r \La &\lra
    M \\
   e_{i} &\mapsto x_{i}  \ea \] is continuous. Since $\bigoplus_{i=1}^r \La$ is
   compact, so is its image $N$.   \epf

\bc \label{fp is compact} Let $M$ be a finitely presented abstract
$\La$-module. Then $M$ is a compact $\La$-module. \ec

\bpf Since $M$ is finitely presented, we have an exact sequence
$\La^{r} \stackrel{f}{\lra} \La^{s} \longrightarrow M
\longrightarrow 0$ for some integers $r$ and $s$. By Lemma
\ref{abstract is cts}(1), the map $f$ is a continuous
$\La$-homomorphism of compact $\La$-modules. Since the category
$\C_{\La}$ is abelian by Lemma \ref{Pontryagin duality}(iv), it
follows that $M$ is an object in $\C_{\La}$.  \epf

\medskip
In view of Corollary \ref{fp is compact}, one may ask the following
two questions. The first is if one can say anything about the
$\mathcal{I}$-adic topology on an abstract $\La$-module $M$. In
general, it is not even clear whether this topology is Hausdorff.
The second question that one may ask is if there are other ways to
endow a finitely presented $\La$-module with a topology such that it
becomes a compact $\La$-module. In response to these two questions,
we have the following proposition. In fact, as we shall see, if $M$
is already a compact $\La$-module, the $\mathcal{I}$-adic topology
is Hausdorff, and it is the only one with which one can endow a
finitely presented $\La$-module in order to make it into a compact
$\La$-module. One may compare the following proposition with
\cite[Prop.\ 5.2.17]{NSW}.

\bp \label{fg compact profinite}
 Let $M$ be a compact $\La$-module. Then
the $\mathcal{I}$-adic topology is finer than the original topology
of $M$, and the canonical homomorphism
\[ \al: M
\lra \plim_{\ai\in\mathcal{I}}M/\ai M\] of $\La$-modules is
injective. Furthermore, if $M$ is a finitely generated $\La$-module,
then the topologies coincide, and the above homomorphism is a
continuous isomorphism of compact $\La$-modules.
 \ep

\bpf: Let $N$ be an open submodule of $M$. Then by continuity, for
each $x\in N$, there exist a neighborhood $V_{x}$ of $x$ and
$\ai_x\in\mathcal{I}$ such that $\ai_x V_{x}\sbs N$. Since $M$ is
compact, it is covered by finitely many such sets, say
$V_{x_1},V_{x_{2}},...,V_{x_{r}}$.  Choose $\ai\in \mathcal{I}$ such
that $\ai \sbs \ai_i$ for all $i=1,...,r$. Then we have $\ai M\sbs
N$, and this shows the first assertion. Since $M$ is Hausdorff under
its original topology, it follows that $M$ is Hausdorff under the
$\mathcal{I}$-adic topology and so
\[\ker\al = \bigcap_{\ai\in\mathcal{I}}\ai M = 0.\]

Now if $M$ is finitely generated, we have a surjection
  \[ \La^{m}\tha
(M~\mathrm{with}~\mathcal{I}\mathrm{-adic}~\mathrm{topology}),\]
which is continuous by Lemma \ref{abstract is cts}(1). This implies
that $M$ with the $\mathcal{I}$-adic topology is compact. By the
first assertion,
 the identity map
\[ (M~\mathrm{with}~\mathcal{I}\mathrm{-}\mathrm{adic}~\mathrm{topology})\lra M\] is continuous. This
in turn gives a continuous bijection between compact spaces and is
therefore a homeomorphism. If $M$ is given the $\mathcal{I}$-adic
topology, then the image of $\al$ is dense in
$\plim_{\ai\in\mathcal{I}}M/\ai M$, and so is surjective since $M$
is compact. \epf

\medskip
We conclude with a description of projective objects in $\C_{\La}$
that are finitely generated over $\La$.

\bp \label{fg projective} Let $P$ be a projective object in
$\C_{\La}$ that is finitely generated over $\La$. Then $P$ is a
projective $\La$-module. Conversely, let $P$ be a finitely generated
projective $\La$-module. Then $P$, endowed with the
$\mathcal{I}$-adic topology, is a compact $\La$-module and is a
projective object in $\C_{\La}$.\ep

\bpf Let $P$ be a projective object in $\C_{\La}$ that is finitely
generated over $\La$. Then there is a surjection $f:\La^{r}\tha P$
of $\La$-modules. By Proposition \ref{fg compact profinite}, the
topology on $P$ is precisely the $\mathcal{I}$-adic topology, and it
follows from Lemma \ref{abstract is cts}(1) that $f$ is a continuous
homomorphism of compact $\La$-modules. Now since $P$ is a projective
object in $\C_{\La}$, the map $f$ has a continuous $\La$-linear
section. In particular, this implies that we have an isomorphism
$\La^r \cong P\oplus(\ker f)$ of $\La$-modules. Hence $P$ is a
projective $\La$-module.

Conversely, suppose that $P$ is a finitely generated projective
$\La$-module. Then there exists a finitely generated projective
$\La$-module $Q$ such that $P\oplus Q$ is a free $\La$-module of
finite rank. We then have a surjection $\pi :\La^n\tha Q$, and this
gives a finite presentation
\[\La^n \lra P\oplus Q \lra P \lra 0\]
of $P$ where the first map sends an element $x$ of $\La^n$ to $(0,
\pi(x))$ and the second map is the canonical projection. It then
follows from Proposition \ref{fp is compact} that $P$ is a compact
$\La$-module under the $\mathcal{I}$-adic topology. Now suppose we
are given the following diagram
 \[ \SelectTips{eu}{} \xymatrix{
                &         P \ar[d]^{\al}    \\
  M  \ar @{->>}[r]^{\e} & N             } \]
of compact $\La$-modules and continuous $\La$-homomorphisms. Since
$P$ is a projective $\La$-module, there is an abstract
$\La$-homomorphism $\be:P \ra M$ such that $\e\be = \al$. On the
other hand, it follows from Lemma \ref{abstract is cts}(3) that
$\be$ is also continuous. Therefore, this shows that $P$ is a
projective object of $\C_{\La}$. \epf

\subsection{Continuous cochains}

\bd Let $G$ be a profinite group. We define $\C_{\La,G}$ to be the
category where the objects are compact $\La$-modules with a
continuous $\La$-linear $G$-action and the morphisms are continuous
$\La[G]$-homomorphisms. Similarly, we define $\Di_{\La,G}$ to be the
category where the objects are discrete $\La$-modules with a
continuous $\La$-linear $G$-action and the morphisms are
(continuous) $\La[G]$-homomorphisms. \ed

\bp \label{Pontryagin duality2} $(i)$ The category $\C_{\La,G}$ is
abelian, has enough projectives and exact inverse limits.

$(ii)$ The category $\Di_{\La,G}$ is abelian, has enough injectives
and exact direct limits.

$(iii)$ The Pontryagin duality induces a contravariant equivalence
between $\C_{\La,G}$ and $\Di_{\La^{\circ}, G}$ $($resp.
$\C_{\La^{\circ},G}$ and $\Di_{\La, G})$. \ep

\bpf We shall prove (iii) first. For a topological group $A$, we
shall denote $A^{\vee}$ to be its Pontryagin dual. By Proposition
\ref{Pontryagin duality}, it suffices to show that if $M$ (resp.,\
$N$) is an object of $\C_{\La,G}$ (resp.,\ $\Di_{\La^{\circ},G}$),
then $M^{\vee}$ (resp.,\ $N^{\vee}$) is an object of
$\Di_{\La^{\circ},G}$ (resp.,\ $\C_{\La,G}$). We define a $G$-action
on $M^{\vee}$ by $\s\cdot f(m) = f(\s^{-1} m)$ for $f\in M^{\vee},
\s\in G$ and $m\in M$. This is clearly $\La^{\circ}$-linear, and
since $G$ is profinite, we may apply \cite[Prop.\ 3]{F} to conclude
that the $G$-action is continuous. The same argument works for $N$.
Hence we have proven (iii). It remains to prove (ii), since (i) will
follow from (ii) and (iii).

To prove (ii), we note that it is clear that $\Di_{\La,G}$ is
abelian and has exact direct limits. It remains to show that it has
enough injectives. By the lemma to follow, we see that the functor
\[ M \mapsto \bigcup_{\ai\in\mathcal{I}}\bigcup_{U}\big(M[\ai]\big)^U :
Mod_{\La[G]}\lra \Di_{\La,G}\] is right adjoint to an exact functor,
and so preserves injectives by \cite[Prop.\ 2.3.10]{Wei}. Since
$Mod_{\La[G]}$ has enough injectives, it follows that $\Di_{\La,G}$
also has enough injectives. \epf

\bl \label{discrete La G} An abstract $\La[G]$-module $N$ is an
object in $\Di_{\La, G}$ if and only if \[ N =
\bigcup_{\ai\in\mathcal{I}}\bigcup_{U}\big(N[\ai]\big)^U,
\]
where $U$ runs through all the open subgroups of $G$. Moreover, if
$M$ is an abstract $\La[G]$-module, then \[
\bigcup_{\ai\in\mathcal{I}}\bigcup_{U}\big(M[\ai]\big)^U,
\] is an object
of $\Di_{\La, G}$, and there is a canonical isomorphism
\[ \Hom_{\La[G],\cts}\Big(N, \bigcup_{\ai\in\mathcal{I}}\bigcup_{U}\big(M[\ai]\big)^U\Big)
\cong \Hom_{\La[G]}(N,M)\] for every $N\in \Di_{\La, G}$. \el

\bpf Suppose $N$ is an object in $\Di_{\La,G}$. Then, in particular,
it is a discrete $\La$-module. By Lemma \ref{discrete La}, we have
$N= \bigcup_{\ai\in\mathcal{I}}N[\ai]$. Let $x\in N[\ai]$. Then by
continuity of the $G$-action, there exists an open subgroup $U$ of
$G$ such that $U\cdot x =x$.

Conversely, suppose that \[ N =
\bigcup_{\ai\in\mathcal{I}}\bigcup_{U}\big(N[\ai]\big)^U.
\] Clearly this
implies that $N= \bigcup_{\ai\in\mathcal{I}}N[\ai]$, and so $N$ is a
discrete $\La$-module. It remains to show that the $G$-action
\[ \theta : G \times N\lra N\]
is continuous. Let $x\in N$, and let $(\s,y)\in\theta^{-1}(x)$. Then
$y\in N[\ai]^U$ for some $\ai\in\mathcal{I}$ and open subgroup $U$.
In particular, we have $(\s,y)\in\s U\times\{y\}\sbs\theta^{-1}(x)$.
Therefore, this proves the first assertion. The second assertion is
an immediate consequence of the first. \epf

\bl \label{compact La G} Let $M$ be an object of $\C_{\La,G}$. Then
$M$ has a fundamental system of neighborhoods of zero consisting of
open $\La[G]$-submodules. \el

\bpf Let $N$ be an open $\La$-submodule of $M$. Then for each $g\in
G$, there exist an open $\La$-submodule $N_{g}$ of $M$ and an open
subgroup $U_{g}$ of $G$ such that $gU_{g}\cdot N_{g}\sbs N$. Since
$G$ is compact, it is covered by finite number of such cosets, say
$g_1U_{g_{1}},...,g_r U_{g_r}$. Set $N_0= \cap_{i=1}^{r}N_{g_i}$.
This is an open $\La$-submodule of $M$. Then $\La[G]\cdot N_0$ is a
$\La[G]$-submodule of $M$ which contains $N_0$ and is contained in
$N$. \epf

\medskip
For the remainder of the subsection, we will be studying the
continuous cochain complex (and its cohomology) of $G$ with
coefficients in certain classes of topological $\La$-modules.

\bd Let $M$ be a topological $\La$-module with a continuous
$\La$-linear $G$-action. The $($inhomogeneous$)$ continuous cochains
$C^{i}_{\cts}(G,M)$ of degree $i\geq 0$ on $G$ with values in $M$
are defined to be the left $\La$-module of continuous maps
$G^{i}\rightarrow M$ with the usual differential
\[(\de^{i} c)(g_{1},..., g_{i+1}) =
g_{1} c(g_{2},...,g_{i+1}) +
\displaystyle\sum_{j=1}^{i}(-1)^{j}c(g_{1},...,g_{j}g_{j+1},...,g_{i+1})
+ (-1)^{i-1}c(g_{1},...,g_{i}),\] which maps
$C^{i}_{\cts}(G,M_{\al})$ to $C^{i+1}_{\cts}(G,M_{\al})$. It then
follows that
\[ \cdots \lra C^{i}_{\cts}(G,M) \stackrel{\de_{M}^{i}}{\lra}
C^{i+1}_{\cts}(G,M)\lra \cdots \] is a complex of $\La$-modules and
its $i$th cohomology group is denoted by $H^{i}_{\cts}(G,M)$. The
following lemma is a standard result (cf.\ \cite[Lemma 2.7.2]{NSW}).
\ed

\bl \label{exact cts cochain} Let \[ 0 \lra M' \stackrel{\al}{\lra}
M \stackrel{\be}{\lra} M'' \lra 0 \] be a short exact sequence of
topological $\La$-modules with a continuous $\La$-linear $G$-action
such that the topology of $M'$ is induced by that of $M$ and such
that $\be$ has a continuous $($not necessarily $\La$-linear$)$
section. Then
\[ 0\ra C^{\bullet}_{\cts}(G,M')\stackrel{\al_{*}}{\ra}
C^{\bullet}_{\cts}(G,M)\stackrel{\be_{*}}{\ra}
C^{\bullet}_{\cts}(G,M'')\ra 0 \] is an exact sequence of complexes
of $\La$-modules. \el

\medskip
We are particularly interested in the case when $M$ is an object of
$\C_{\La,G}$ or $\Di_{\La,G}$. We now discuss cohomology and limits.

\bp \label{direct limit cochain} Let $N = \ilim_{\al} N_{\al}$ be an
object of $\Di_{\La,G}$, where $N_{\al}\in\Di_{\La,G}$. Then we have
an isomorphism
\[ C^i_{\cts}(G,N) \cong \ilim_{\al}C^i_{\cts}(G,N_{\al}) \] of continuous
cochain groups which induces an isomorphism
\[ H^i_{\cts}(G,N)\cong \ilim_{\al}H^i_{\cts}(G,N_{\al})\] of
cohomology groups. \ep

\bpf The first isomorphism is immediate and the second follows from
the first since direct limit is exact. \epf

\medskip
In the next proposition, we shall examine the relationship between
cohomology and inverse limit. We shall denote ${\plim}^{(i)}$ to be
the $i$th derived functor of $\plim$.

\bp \label{inverse limit cochain} Let $M = \plim_{\al} M_{\al}$ be
an object in $\C_{\La,G}$, where each $M_{\al}$ is finite. Then we
have an isomorphism
\[ C_{\cts}(G,M) \cong \plim_{\al}C_{\cts}(G,M_{\al}) \] of complexes of $\La$-modules and a spectral sequence
\[ {\plim_{\al}}^{(i)}H^j_{\cts}(G, M_{\al}) \Longrightarrow H^{i+j}_{\cts}(G, M). \]
Suppose further that $G$ has the property that $H^{m}_{\cts}(G,N)$
is finite for all finite discrete $\La$-modules $N$ with a
continuous commuting $G$-action and for all $m\geq 0$. Then
\[H^{i}_{\cts}(G,M) \cong \plim_{\al}H^{i}_{\cts}(G,M_{\al}).\] \ep

\bpf  The first assertion is immediate from the definition. The
second assertion follows from a similar argument as in \cite[Prop.\
8.3.5]{Ne}. We consider the two hypercohomology spectral sequences
for the functor $\plim$ and the inverse system $C^i_{\cts}(G,
M_{\al})$:
\[ \ba{c} {\plim_{\al}}^{(j)} C^i_{\cts}(G, M_{\al}) \Longrightarrow
H^{i+j} \\
 {\plim_{\al}}^{(i)} H^j_{\cts}(G, M_{\al}) \Longrightarrow
H^{i+j}.  \ea \] For each $i$, it is clear that
\[ \plim_{\al}C^i_{\cts}(G, M_{\al}) \lra C^i_{\cts}(G, M_{\al})\] is surjective
for every $\al$, and so the inverse system $C^i_{\cts}(G,M_{\al})$
is ``weakly flabby" in the sense of \cite[Lemma 1.3]{Je}. Therefore,
by \cite[Thm.\ 1.8]{Je}, we have that ${\plim_{\al}}^{(j)}
C^i_{\cts}(G, M_{\al}) = 0$ for $j>0$. Hence, the first spectral
sequence degenerates and we obtain
 \[H^i = H^i\big(\plim_{\al} C_{\cts}(G, M_{\al})\big) \cong H^i\big( C_{\cts}(G,
 M)\big) = H^i_{\cts}(G,M).\]

For the last assertion, the additional assumption allows one to
invoke \cite[Cor.\ 7.2]{Je} to conclude that
${\plim_{\al}}^{(i)}H^{j}_{\cts}(G, M_{\al}) =0$ for $i>0$. \epf

 For the remainder of the subsection, we let $\mathcal{A}$ denote
either $\C_{\La,G}$ or $\Di_{\La,G}$. Let $M^{\bullet}$ be a complex
of objects in $\mathcal{A}$ with differentials denoted by
$d_{M}^{i}$. We define $C^{\bullet}_{\cts}(G,M^{\bullet})$ by
\[ C^{n}_{\cts}(G,M^{\bullet}) =
\bigoplus_{i+j=n}C^{j}_{\cts}(G,M^{i}) .\] Its differential
$\de_{M^{\bullet}}^{i+j}$ is determined as follows: restriction of
$\de_{M^{\bullet}}^{i+j}$ to $C^{j}_{\cts}(G,M^{i})$ is the sum of
\[ (d_{M}^{i})_{*} : C^{j}_{\cts}(G,M^{i})\lra
C^{j}_{\cts}(G,M^{i+1}) \] and \[ (-1)^{i}\de_{M^{i}}^{j} :
C^{j}_{\cts}(G,M^{i}) \lra C^{j+1}_{\cts}(G,M^{i}). \] We denote its
$i$th cohomology group by $H^{i}(G,M^{\bullet})$.

\bp \label{exact cochain compact discrete} Let $0\ra
M'\stackrel{\al}{\ra} M \stackrel{\be}{\ra} M''\ra 0$ be an exact
sequence of objects in $\mathcal{A}$. Then
\[ 0\ra C^{\bullet}_{\cts}(G,M')\stackrel{\al_{*}}{\ra}
C^{\bullet}_{\cts}(G,M)\stackrel{\be_{*}}{\ra}
C^{\bullet}_{\cts}(G,M'')\ra 0 \] is an exact sequence of complexes
of $\La$-modules. The statement also holds true if we replace
$M',M,M''$ by complexes of objects in $\mathcal{A}$.\ep

\bpf By Lemma \ref{exact cts cochain}, it suffices to show that
$\be$ has a continuous section. If $\mathcal{A} = \Di_{\La,G}$, this
is obvious. In the case when $\mathcal{A} = \C_{\La,G}$, since every
compact $\La$-module is profinite by Proposition \ref{Pontryagin
duality}, every continuous surjection has a continuous section. \epf

Let $M^{\bullet}$ be a complex of objects in $\mathcal{A}$. The
filtration $\tau_{\leq j}M^{\bullet}$ induces a filtration
 \[ \tau_{\leq j}C^{\bullet}_{\cts}(G,M^{\bullet}) = C^{\bullet}_{\cts}(G,\tau_{\leq j}
 M^{\bullet})\] on the cochain groups which fit into the following
 exact sequence of complexes
 \[ 0\lra C^{\bullet}_{\cts}(G,\tau_{\leq j}M^{\bullet}) \lra
 C^{\bullet}_{\cts}(G,\tau_{\leq j+1}M^{\bullet}) \lra
 \tau_{\leq j+1}C^{\bullet}_{\cts}(G,M^{\bullet})/\tau_{\leq j}C^{\bullet}_{\cts}(G,M^{\bullet})
 \lra 0\]
by Proposition \ref{exact cochain compact discrete}. This filtration
gives rise to the following hypercohomology spectral sequence
\[H^{i}_{\cts}\big(G,H^{j}(M^{\bullet})\big)\Longrightarrow
H^{i+j}_{\cts}(G,M^{\bullet}),\] which is convergent if
$M^{\bullet}$ is cohomologically bounded below.

\bl \label{bounded below qis} Let $f:M^{\bullet}\lra N^{\bullet}$ be
a quasi-isomorphism of cohomologically bounded below complexes of
objects in $\mathcal{A}$. Then the induced map
\[ f_{*}:C^{\bullet}_{\cts}(G,M^{\bullet})\lra
C^{\bullet}_{\cts}(G,N^{\bullet}) \] is also a quasi-isomorphism.
\el

\bpf The map $f$ induces isomorphisms
\[ H^{i}_{\cts}\big(G,H^{j}(M^{\bullet})\big) \stackrel{\sim}{\lra}
H^{i}_{\cts}\big(G,H^{j}(N^{\bullet})\big). \] By convergence of the
above spectral sequence, this implies that the induced maps
\[H^{i}_{\cts}(G,M^{\bullet})\lra H^{i}_{\cts}(G,N^{\bullet})\] are
isomorphisms. \epf

Hence we can conclude the following.

\bp \label{derived cts cochain} The functor
\[ C_{\cts}^{\bullet}(G,-):\mathrm{Ch}^+(\mathcal{A})\lra
\mathrm{Ch}^+(Mod_{\La})\] preserves homotopy, exact sequences and
quasi-isomorphisms, hence induces the following exact derived
functors
\[ \ba{c}
 \R\Ga_{\cts}(G,-): \DD^{b}(\C_{\La,G})\longrightarrow
\DD^{+}(Mod_{\La})\\
 \R\Ga_{\cts}(G,-): \DD^{+}(\Di_{\La,G})\longrightarrow
 \DD^{+}(Mod_{\La}). \ea\] \ep

\bpf This proposition follows from what we have done so far. The
only subtlety lies in the fact that $\C_{\La,G}$ does not
necessarily have enough injectives and therefore we do not know if
$\DD^{+}(\C_{\La,G})$ exists. However, we know that $\C_{\La,G}$ has
enough projectives. Therefore, $\DD^-(\C_{\La,G})$ exists, and we
may apply Lemma \ref{bounded below qis} to $\DD^b(\C_{\La,G})$. \epf

We now like to extend Proposition \ref{inverse limit cochain} to the
case of complexes. Before that, we first prove a lemma which will be
required in our discussion.

\bl Let $f:M\lra N$ be a morphism of objects in $\C_{\La,G}$. Then
there exists a directed indexing set $I$ with the following
properties$:$

$(1)$ There exist a fundamental system  $\{U_{i}\}$ $($resp.,
$\{V_i\})$ of neighborhoods of zero consisting of open
$\La[G]$-submodules of $M$ $($resp., $N)$.

$(2)$ For each $i\in I$, there is a $\La[G]$-homomorphism $f_i :
M/U_i\lra N/V_i$ which fits into the following commutative diagram
 \[  \SelectTips{eu}{} \xymatrix{
   M \ar[d]_{} \ar[r]^{f} & N \ar[d]^{} \\
   M/U_i \ar[r]^{f_i} & N/V_i   }\]
where the vertical morphisms are the canonical quotient map.

$(3)$ One has  $f = \plim_i f_i$. \el

\bpf Let $\{U_{\al}\}_{\al\in I_M}$ (resp., $\{V_{\be}\}_{\be\in
I_N}$) be a system of neighborhoods of zero consisting of open
$\La[G]$-submodules of $M$ $($resp., $N)$. Then we set $I= I_M\times
I_N, U_{\al,\be} = U_{\al}\cap f^{-1}(V_{\be})$ and $V_{\al,\be} =
V_{\be}$. It is then straightforward to verify that $f$ factors
through $M/U_{\al,\be}$ to give a $\La[G]$-homomorphism $f_{\al,\be}
: M/U_{\al,\be}\lra N/V_{\al,\be}$ and $f =
\plim_{\al,\be}f_{\al,\be}$. \epf

In view of the above lemma, we say that a morphism $f :M\lra N$ in
$\C_{\La,G}$ is compatible with a directed indexing set $I$ if the
conclusion in the lemma holds. By the lemma, we have that for every
morphism $f:M\lra N$ in $\C_{\La,G}$, there exists a directed
indexing set $I$ such that $f$ is compatible with $I$. In
particular, if $M$ is a bounded complex in $\C_{\La,G}$, we can find
a directed indexing set $I$ such that the differentials are
compatible with $I$.

\bp \label{inverse limit cochain2} Suppose that $G$ has the property
that $H^{m}_{\cts}(G,N)$ is finite for all finite discrete
$\La$-modules $N$ with a continuous commuting $G$-action and for all
$m\geq 0$. Let $M^{\bullet} = \plim_{i\in I} M^{\bullet}_i$ be a
bounded complex of objects in $\C_{\La,G}$ with $I$-compatible
differentials. Then we have the following isomorphism
\[H^{n}_{\cts}(G,M^{\bullet}) \cong \plim_{i}H^{n}_{\cts}(G,M^{\bullet}_i)\]
of hypercohomology groups for each $n$. \ep

\bpf The canonical chain map $M^{\bullet}\lra M^{\bullet}_i$ induces
the following morphism of (convergent) spectral sequences
\[  \SelectTips{eu}{} \xymatrix{
  H^{r}_{\cts}\big(G,H^{s}(M^{\bullet})\big)\Longrightarrow
H^{r+s}_{\cts}(G,M^{\bullet}) \ar[d]^{} \\
  H^{r}_{\cts}\big(G,H^{s}(M_i^{\bullet})\big)\Longrightarrow
H^{r+s}_{\cts}(G,M_i^{\bullet}) }\] which is compatible with $i$. By
the hypothesis, the bottom spectral sequence is a spectral sequence
of finite $\La$-modules. Therefore, the inverse limit is compatible
with the inverse system of the spectral sequences, and we have the
following morphism
\[ \SelectTips{eu}{} \xymatrix{
H^{r}_{\cts}\big(G,H^{s}(M^{\bullet})\big)\Longrightarrow
H^{r+s}_{\cts}(G,M^{\bullet}) \ar[d]^{} \\
 \plim_i H^{r}_{\cts}\big(G,H^{s}(M_i^{\bullet})\big)\Longrightarrow
\plim_i H^{r+s}_{\cts}(G,M_i^{\bullet}) }\] of (convergent) spectral
sequences. By Proposition \ref{inverse limit cochain}, we have the
isomorphisms
 \[ H^{r}_{\cts}\big(G,H^{s}(M^{\bullet})\big) \cong \plim_i H^{r}_{\cts}\big(G,H^{s}(M_i^{\bullet})\big).\]
Hence, by the convergence of the spectral sequences, we obtain the
required isomorphism. \epf

\medskip
For ease of notation, we will drop the `$\bullet$' for complexes. We
also drop the notation `cts'. Therefore, we write $C(G, M)$ as the
complex of continuous cochains and $\R\Ga(G,M)$ for its derived
functor. Its $i$th cohomology group is then written as $H^{i}(G,M)$.

\subsection{Total cup products}  \label{cup product section}

We first recall the definition for topological $G$-modules (in other
words, abelian Hausdorff topological groups with a continuous
$G$-action).

\bd (Cup products) Let $A, B$ and $C$ be topological $G$-modules.
Suppose
\[ \langle\ ,\ \rangle : A\times B \lra C \] is a continuous map satisfying
$\s\langle a,b\rangle = \langle\s a, \s b\rangle$ for $a\in A, b\in
B$ and $\s\in G$. Then we define the cup product on the cochain
groups
\[ C^i(G,A)\times C^j(G,B)\lra C^{i+j}(G,C)\]
as follows: for $\al\in C^i(G,A), \be\in C^j(G,B)$ and
$\s_{1},...,\s_{i+j}\in G$, we have
\[(\al\cup \be)(\s_{1},...,\s_{i+j}) = \Big\langle\al(\s_{1},...,\s_{i}),
\s_{1}\cdots\s_{i}\be(\s_{i+1},...,\s_{i+j})\Big\rangle.\]
 The cup product satisfies the following relation
 \[ \de_{C}(\al\cup\be) = (\de_{A}\al)\cup\be + (-1)^i\al\cup(\de_B\be)\]
and induces a pairing
\[  H^i(G,A)\times H^j(G,B)\lra H^{i+j}(G,C)\]
on the cohomology groups. \ed

\medskip
Now fix a prime $p$. For the remainder of the paper, we shall assume
that our profinite ring $\La$ is pro-$p$. In other words, for each
$\ai\in\mathcal{I}$, the ring $\La/\ai$ is finite of a $p$-power
cardinality. Let $M$ and $N$ be objects in $\C_{\La,G}$ and
$\Di_{\La^{\circ},G}$ respectively, and let $A$ be a topological
$G$-module. Suppose there is a continuous pairing
\[ \langle\ ,\ \rangle : N\times M \lra A
\] such that

(1) $\s\langle y,x\rangle = \langle\s y, \s x\rangle$ for $x\in M,
y\in N$ and $\s\in G$, and

(2) $\langle y\la, x\rangle = \langle y,\la x\rangle$ for $x\in M,
y\in N$ and $\la\in \La$.

As before, condition (1) will give rise to the cup product
\[ C^i(G,N)\times C^j(G,M)\lra C^{i+j}(G,A),\]
which is $\La$-balanced by condition (2). The cup product induces a
group homomorphism
\[ C^i(G,N)\ot_{\La} C^j(G,M)\lra C^{i+j}(G,A)\]
which gives rise to the following morphism
\[ C(G,N)\ot_{\La} C(G,M)\lra C(G,A)\] of complexes of abelian
groups. Taking the adjoint, we have a morphism
\[ C(G,M) \lra \Hom_{\Zp}\big(C(G,N), C(G,A)\big) \] of complexes of
$\La$-modules.

\bl \label{cup product1} Suppose we are given another continuous
pairing
\[ (\ ,\ ): N'\times M' \lra A
\] such that
$(1)$ $\s( y',x') = (\s y', \s x')$ for $x'\in M', y'\in N'$ and
$\s\in G$;

$(2)$ $( y'\la, x') = ( y',\la x')$ for $x'\in M', y'\in N'$ and
$\la\in \La$, and

$(3)$ there are morphisms $f:N'\lra N$ in $\Di_{\La^{\circ},G}$ and
$g:M\lra M'$ in $\C_{\La,G}$ such that the following diagram
\[  \SelectTips{eu}{} \xymatrix @C=0.7in{
  N'\ot_{\La}M \ar[d]_{f\ot\mathrm{id}} \ar[r]^{\mathrm{id}\ot g} & N'\ot_{\La}M' \ar[d]^{(~,~)} \\
  N\ot_{\La}M \ar[r]^{\langle~,~\rangle} & A   } \] commutes. Then
we have the following commutative diagram
\[  \SelectTips{eu}{} \xymatrix{
  C(G,M) \ar[d]_{g_*} \ar[r]^{} & \Hom_{\Zp}\big(C(G,N), C(G,A)\big) \ar[d]^{f_*} \\
  C(G,M') \ar[r]^{} & \Hom_{\Zp}\big(C(G,N'), C(G,A)\big)   } \] of complexes
of $\La$-modules. \el

\bpf It follows from a direct calculation that following diagram
\[  \SelectTips{eu}{} \xymatrix @C=0.7in{
  C(G,N')\ot_{\La}C(G, M) \ar[d]_{f\ot\mathrm{id}} \ar[r]^{\mathrm{id}\ot g} &
  C(G, N')\ot_{\La} C(G, M') \ar[d]^{\cup_{(~,~)}} \\
  C(G, N)\ot_{\La}C(G, M) \ar[r]^{\cup_{\langle~,~\rangle}} & C(G, A)   } \]
is commutative, where $\cup_{(~,~)}$ and $\cup_{\langle~,~\rangle}$
are the cup products induced by the pairings $(~,~)$ and
$\langle~,~\rangle$ respectively. By taking the adjoint and another
straightforward calculation, we have the commutative diagram in the
lemma. \epf

Now let $M$ and $N$ be bounded complexes of objects in $\C_{\La,G}$
and $\Di_{\La^{\circ},G}$ respectively, and let $A$ be a bounded
complex of topological $G$-modules. Suppose there is a collection of
continuous pairings
\[ \langle\ ,\ \rangle_{a,b} : N^{a}\times M^b\lra A^{a+b}\]
where each pairing satisfies conditions (1) and (2), and the
following hold:

(a) $\langle d^a_N y,x\rangle_{a+1,b} = d^{a+b}_A\big(\langle y,
x\rangle_{a,b}\big)$ for $y\in N^a$ and $x\in M^b$, and

(b) $(-1)^a\langle y, d^b_Mx\rangle_{a,b+1} = d^{a+b}_A\big(\langle
y, x\rangle_{a,b}\big)$ for $y\in N^a$ and $x\in M^b$.

\noindent For each pair $(a,b)$, we have a morphism
\[ \cup_{ij}^{ab} : C^{i}(G,N^a) \ot_{\La} C^{j}(G,M^b) \lra
C^{i+j}(G, A^{a+b}) \] of abelian groups induced by the cup product.
Then the total cup product
\[ \cup : C(G,N)\ot_{\La} C(G,M) \lra C(G,A)\] is a morphism of complexes of $\Zp$-modules
given by the collection $\cup = \big( (-1)^{ib}\cup^{ab}_{ij}\big)$.
The definition given for the total cup products follows that in
\cite[3.4.5.2]{Ne}. We also have an analogous result to Lemma
\ref{cup product1} for complexes.

\bl \label{cup product2} Suppose we are given another collection of
continuous pairings
\[ \langle\ ,\ \rangle_{a,b} : N'^{a}\times M'^b\lra A'^{a+b}\]
as above. Then we have the following commutative diagram
\[ \SelectTips{eu}{} \xymatrix{
  C(G,M) \ar[d]_{g_*} \ar[r]^{} & \Hom_{\Zp}\big(C(G,N), C(G,A)\big) \ar[d]^{f_*} \\
  C(G,M') \ar[r]^{} & \Hom_{\Zp}\big(C(G,N'), C(G,A)\big)   } \] of complexes
of $\La$-modules. \el

\subsection{Tate cohomology groups} \label{Tate cohomology groups}

We shall now describe the Tate cochain complexes of a finite group
$G$. We begin by giving an alternative description of the
(inhomogeneous) cochain complexes. Throughout this subsection, $G$
will always denote a finite group. Consider the standard
$\Z[G]$-resolution (in inhomogeneous form) of $\Z$ (cf.\
\cite[Sect.\ 6.5]{Wei})\footnote{Weibel calls this the unnormalized
bar resolution.},
\[ X_0 \longleftarrow X_1 \longleftarrow X_2
\longleftarrow \cdots\] where $X_0 = \Z[G]$ and, for $n\geq 1$,
$X_n$ is the free $\Z[G]$-module generated by the set of all symbols
$(g_1,...,g_n)$ with $g_i\in G$, and the differentials are given by
the formula
\[ \partial^n(g_1,...,g_n) = g_{1}(g_{2},...,g_{n}) +
\sum_{j=1}^{n-1}(-1)^{j}(g_{1},...,g_{j}g_{j+1},...,g_{n}) +
(-1)^{n}(g_{1},...,g_{n-1}). \]

\noindent For any $\Z[G]$-module $M$, there is a natural isomorphism
$C^i(G,M)\lra \Hom_{\Z[G]}(X_i, M)$ which is compatible with the
differentials, thus giving an identification of complexes.
Furthermore, if $M$ is a $\La[G]$-module, the above identification
is an isomorphism of complexes of $\La$-modules.

We now construct the complete cochain groups. For a $\Z[G]$-module
$A$, we write $A^* = \Hom_{\Z}(A,\Z)$. Note that this is a
$\Z[G]$-module in a natural way. Applying $\Hom_{\Z}(-,\Z)$ to the
long exact sequence
\[ 0 \longleftarrow \Z\longleftarrow  X_0 \longleftarrow X_1 \longleftarrow X_2
\longleftarrow \cdots,\]  we obtain the following long exact
sequence
\[ 0 \longrightarrow \Z \lra X_0^* \longrightarrow X_1^* \longrightarrow X_2^*
\longrightarrow \cdots,\] since each $X_i$ is a free $\Z$-module.
Splicing the two long exact sequence and applying $\Hom_{\Z[G]}(-,
M)$ to the resulting long exact sequence, we obtain the following
complex
\[ \cdots \lra \Hom_{\Z[G]}(X_1^*,
M) \lra \Hom_{\Z[G]}(X_0^*, M) \lra \Hom_{\Z[G]}(X_0, M)\lra \cdots
\]

The completed cochain complexes $\hat{C}^i(G,M)$ are defined by
\begin{equation*}  \hat{C}^i(G,M)=
\begin{cases} C^i(G,M) \cong \Hom_{\Z[G]}(X_i, M) & \text{\mbox{if} $i\geq 0$},\\ \Hom_{\Z[G]}(X_{-1-i}^*, M) & \text{\mbox{if} $i\leq -1$.}
\end{cases} \end{equation*}

Following \cite[5.7.2]{Ne}, we may extend the above definition to a
complex $M^{\bullet}$ of $G$-modules by setting
\[ \hat{C}^n(G,M^{\bullet}) = \bigoplus_{i+j = n}\hat{C}^i(G, M^j)\]
with differential defined using the sign conventions of the previous
sections. As before, for ease of notation, we will drop the
`$\bullet$' for complexes. The usual cup product for Tate cohomology
groups (cf. \cite[Prop.\ 1.4.6]{NSW}) extends to a total cup product
with the same sign convention as in the preceding section.

\section{Duality over pro-$p$ rings} \label{duality}

Let $p$ be a fixed prime. Throughout the section, our profinite ring
$\La$ will always be pro-$p$. In this section, we will formulate and
prove Tate's (and Poitou's) local and global duality theorems.

\subsection{Tate's local duality}

Let $F$ be a nonarchimedean local field with characteristic not
equal to $p$. Fix a separable closure $F^{\mathrm{sep}}$ of $F$. Set
$G_{F} = \mathrm{Gal}(F^{\mathrm{sep}}/F)$.

\bl \label{local class field theory} We have
\begin{equation*}  H^{j}(G_{F},\Qp/\Zp(1))\cong
\begin{cases} \Qp/\Zp & \text{\mbox{if} $j=2$},\\  0 & \text{\mbox{if} $j>2$.}
\end{cases} \end{equation*} \el

\bpf For $j>2$, the conclusion follows from the fact that $G_{F}$
has $p$-cohomological dimension 2 (see \cite[Thm.\ 7.1.8(i)]{NSW}).
By \cite[Thm.\ 7.1.8(ii)]{NSW}, we have
$H^{2}(G_{v},\Z/p^{r}(1))\cong \Z/p^{r}$. The assertion now follows
by taking direct limits. \epf

By the preceding lemma, we have a quasi-isomorphism $\Qp/\Zp[-2]
\stackrel{i}{\lra}\tau_{\geq 2}C(G_{F}, \Qp/\Zp(1))$ of complexes of
$\Zp$-modules. Since $\Qp/\Zp$ is an injective $\Zp$-module, the map
$i$ has a homotopy inverse. We shall fix one such map
\[r:\tau_{\geq 2}C(G_{F}, \Qp/\Zp(1))\lra \Qp/\Zp[-2]. \] This
gives a morphism
\[ \theta: C(G_F,\Qp/\Zp(1))\lra \tau_{\geq
2}C(G_{F}, \Qp/\Zp(1))\stackrel{r}{\lra} \Qp/\Zp[-2]\] of complexes
of $\Zp$-modules.

\medskip Let $M$ be a bounded complex of objects in $\C_{\La,G_F}$. We shall write $M^{\vee}$
to be the complex $\Hom_{\cts}(M, \Qp/\Zp)$. The obvious pairing
\[ M^{\vee}(1)\ot_{\La} M \lra \Qp/\Zp\] induces the total cup
product
\[ C\big(G_F, M^{\vee}(1)\big) \ot_{\La} C(G_F, M) \lra C\big(G_F,
\Qp/\Zp(1)\big). \] Suppose that $N$ is another bounded complex of
objects in $\C_{\La,G_F}$, and there is a morphism $f: M\lra N$ of
complexes in $\C_{\La,G_F}$. Then we have the following commutative
diagram
\[ \SelectTips{eu}{} \xymatrix @C=0.7in{
  N^{\vee}(1)\ot_{\La}M \ar[d]_{f^{\vee}\ot\mathrm{id}} \ar[r]^-{\mathrm{id}\ot f} &
  N^{\vee}(1)\ot_{\La}N \ar[d] \\
  M^{\vee}(1)\ot_{\La}M \ar[r] & \Qp/\Zp(1)  } \]
with the obvious pairings. Applying cochains and $\theta$, we obtain
the following commutative diagram
\[  \SelectTips{eu}{} \xymatrix @C=0.6in{
  C\big(G_F, N^{\vee}(1)\big)\ot_{\La}C(G_{F},M) \ar[d]_{f^{\vee}\ot\mathrm{id}} \ar[r]^-{\mathrm{id}\ot f} &
  C\big(G_F, N^{\vee}(1)\big)\ot_{\La}C(G_F, N) \ar[d] \\
  C\big(G_{F},M^{\vee}(1)\big)\ot_{\La}C(G_F, M) \ar[r] & \Qp/\Zp[-2]  } \]
which induces the following commutative diagram
\[  \SelectTips{eu}{} \xymatrix @C=0.6in  @R=0.3in{
  C(G_F, M) \ar[d] \ar[r]^-{\al_M} & \Hom_{\Zp}\Big(C\big(G_F,M^{\vee}(1)\big),\Qp/\Zp\Big)[-2] \ar[d] \\
  C(G_F, N) \ar[r]^-{\al_N} & \Hom_{\Zp}\Big(C\big(G_F, N^{\vee}(1)\big),\Qp/\Zp\Big)[-2]  } \]
of complexes of $\La$-modules by Lemma \ref{cup product1}. We are
now able to prove the following formulation of Tate's local duality.

\bt \label{Tate local adic} Let $M$ be a bounded complex of objects
in $\C_{\La,G_F}$. Then we have the following isomorphism
\[ \R\Ga(G_{F},M)\lra
\R\Hom_{\Zp}\Big(\R\Ga\big(G_{F},M^{\vee}(1)\big),\Qp/\Zp\Big)[-2]
\] in $\DD(Mod_{\La})$. \et

\bpf  We shall show that $\al_M$ (in the above diagram) is a
quasi-isomorphism. Now if $A\lra B\lra C\lra A[1]$ is an exact
triangle \footnote{
    We write an exact triangle
    $A \to B \to C \to A[1]$ more compactly as $A \to B \to C$ throughout.}
     in $\DD^b(\C_{\La,G_F})$, we then have a morphism
\[ \SelectTips{eu}{} \entrymodifiers = {!! <0pt,.5ex>+} \xymatrix @C=0.6in{
 \R\Ga(G_{F},A)
  \ar[d]_{} \ar[r]^-{\al_A} & \R\Hom_{\Zp}\Big(\R\Ga\big(G_{F},A^{\vee}(1)\big),\Qp/\Zp\Big)[-2]
\ar[d]^{} \\
\R\Ga(G_{F},B) \ar[r]^-{\al_B}\ar[d] &
\R\Hom_{\Zp}\Big(\R\Ga\big(G_{F},B^{\vee}(1)\big),\Qp/\Zp\Big)[-2] \ar[d]^{} \\
 \R\Ga(G_{F},C) \ar[r]^-{\al_{C}}
& \R\Hom_{\Zp}\Big(\R\Ga\big(G_{F},C^{\vee}(1)\big),\Qp/\Zp\Big)[-2]
} \] of exact triangles. Therefore, if any two of the morphisms
$\al_{A}, \al_B$ and $\al_C$ are isomorphisms, so is the third. For
a bounded complex $M$ in $\C_{\La,G_F}$, we have the following exact
triangle \[ \s_{\leq i-1}M \lra \s_{\leq i}M \lra M^i[-i] .
\] Therefore, by induction, we are reduced to showing
that $\al_M$ is a quasi-isomorphism in the case when $M$ is a single
module. Write $M = \plim_{\be} M_{\be}$, where each $M_{\be}$ is a
finite module. By the functoriality of $\al$, we have the following
commutative diagram
\[  \SelectTips{eu}{} \entrymodifiers={!! <0pt, .8ex>+}  \xymatrix @C=0.8in{
  C(G_F, M) \ar[d]_{u} \ar[r]^-{\al_M} & \Hom_{\Zp}\Big(C\big(G_F,M^{\vee}(1)\big),\Qp/\Zp\Big)[-2] \ar[d]^{v} \\
  \plim_{\be} C(G_F, M_{\be}) \ar[r]^-{\plim\al_{M_{\be}}} &
  \plim_{\be} \Hom_{\Zp}\Big(C\big(G_F,M_{\be}^{\vee}(1)\big),\Qp/\Zp\Big)[-2]  } \]
of complexes of $\La$-modules. By Proposition \ref{direct limit
cochain} and Proposition \ref{inverse limit cochain}, we have that
$u$ and $v$ in the above diagram are isomorphisms of complexes, and
the vertical maps in the following commutative diagram
\begin{equation*}  \SelectTips{eu}{} \entrymodifiers={!! <0pt, .8ex>+}  \xymatrix @C=1.1in{
  H^i(G_F, M) \ar[d]_{u_*} \ar[r]^-{(\al_M)_*} &
  \Hom_{\Zp}\Big(H^{2-i}\big(G_F,M^{\vee}(1)\big),\Qp/\Zp\Big) \ar[d]^{v_*} \\
  \plim_{\be} H^i(G_F,M_{\be}) \ar[r]^-{\plim(\al_{M_{\be}})_*} &
   \plim_{\be} \Hom_{\Zp}\Big(H^{2-i}\big(G_F, M_{\be}^{\vee}(1)\big),\Qp/\Zp\Big)  }
   \end{equation*} are isomorphisms. Since each $(\al_{M_{\be}})_*$ is an isomorphism by Tate
local duality \cite[Thm.\ 7.2.6]{NSW}, we have the required
conclusion. \epf

\subsection{Global duality over pro-$p$ rings}

Let $F$ be a global field with characteristic not equal to $p$, and
let $S$ be a finite set of primes of $F$ containing all primes above
$p$ and all archimedean primes of $F$ (if $F$ is a number field).
Let $S_{f}$ (resp.,\ $S_{\mathds{R}}$) denote the collection of
non-archimedean primes (resp.,\ real primes) of $F$ in $S$.

Fix a separable closure $F^{\mathrm{sep}}$ of $F$. Set $G_{F,S} =
\mathrm{Gal}(F_{S}/F)$, where $F_{S}$ is the maximal subextension of
$F^{\mathrm{sep}}/F$ unramified outside $S$. For each $v\in S_{f}$,
we fix a separable closure $F_{v}^\mathrm{sep}$ of $F_{v}$ and an
embedding $F^{\mathrm{sep}} \hookrightarrow F_{v}^\mathrm{sep}$.
This induces a continuous group homomorphism
$G_{v}:=\mathrm{Gal}(F_{v}^\mathrm{sep}/F_{v})\ra G_{F,S}$. If $v$
is a real prime, we also write $G_v$ for
$\mathrm{Gal}(\mathds{C}/\mathds{R})$.

If $M$ is a complex in $\C_{\La,G_{F,S}}$ (resp.,
$\Di_{\La,G_{F,S}}$), then we can view $M$ as a complex in
$\C_{\La,G_v}$ (resp., $\Di_{\La,G_v}$) via the continuous
homomorphism $G_{v}\lra G_{F,S}$. Therefore, the cochain complexes
$C(G_{F,S},M)$ and $C(G_{v},M)$ can be defined. Recall that for
$v\in S_f$, we have the restriction map
\[\mathrm{res}_v : C(G_{F,S},M)\lra C(G_v,M)
\] induced by the group homomorphism $G_v\lra G_{F,S}$. For a real
prime $v$, we have the following
\[ \mathrm{res}_v : C(G_{F,S},M) \lra C(G_v, M) \hookrightarrow \hat{C}(G_v,
M).\] To shorten notation in what follows, for $v \in
S_{\mathds{R}}$, we will abuse notation and use $C(G_v,M)$,
$H^i(G_v,M)$, and $\R\Ga(G_v,M)$ to denote the Tate cochains
$\hat{C}(G_v,M)$, its cohomology groups, and its derived object. We
now make the following definition.

\bd \label{def compact support} Let $M$ be a complex in
$\C_{\La,G_{F,S}}$ or $\Di_{\La,G_{F,S}}$. The complex of continuous
cochains of $M$ with compact support is defined as
\[ C_{c}(G_{F,S},M) =
\Cone\left(C(G_{F,S},M)\stackrel{\mathrm{res}_S}{\lra}\bigoplus_{v\in
S}C(G_{v},M)\right)[-1]~, \] where the elements of
\[C^{i}_{c}(G_{F,S},M) = C^{i}(G_{F,S},M)\oplus \left( \bigoplus_{v\in
S}C^{i-1}(G_{v},M)\right)
\] have the form $(a,a_{S})$ with $a\in C^{i}(G_{F,S},M),
a_{S}=(a_{v})_{v\in S}, a_v\in C^{i-1}(G_{v},M)$, and the
differential is given by
\[ d(a,a_{S}) = (da, -\mathrm{res}_{S}(a) - da_{S}).\]
The $i$th cohomology group of $C_{c}(G_{F,S},M)$ is denoted by
$H_{c}^{i}(G_{F,S},M)$.\ed

\br If $F$ is a function field in one variable over a finite field
or $F$ is a totally imaginary number field, then $S_{\mathbb{R}}$ is
empty, and the cone is given by\[
\Cone\left(C(G_{F,S},M)\stackrel{\mathrm{res}_{S_f}}{\lra}\bigoplus_{v\in
S_{f}}C(G_{v},M)\right)[-1].\] Now suppose that $p$ is odd and $F$
is a number field with at least one real prime. Let $v\in
S_{\mathbb{R}}$. Then $\hat{H}^i(G_v, M) = 0$ for every $M$ in
$\C_{\La,G_{F,S}}$ (resp., $\Di_{\La,G_{F,S}}$) and for all $i$,
since $G_v$ is a finite group of order 2 and $M$ is an inverse limit
of finite $p$-groups (resp., direct limit of finite $p$-groups).
Therefore, it follows that the canonical map
\[ \Cone\left(C(G_{F,S},M)\stackrel{\mathrm{res}_{S_f}}{\lra}\bigoplus_{v\in
S_{f}}C(G_{v},M)\right)[-1]\lra C_{c}(G_{F,S}, M) \] is a
quasi-isomorphism. Therefore, one may take the above cone as a
definition of the complex of continuous cochains with compact
support in this case.   \er

\bp \label{derived cts compact cochain} The functor
\[ \ba{c} C_c(G_{F,S},-):\mathrm{Ch}^+ (\C_{\La, G_{F,S}})\lra
\mathrm{Ch}(Mod_{\La}) \\
 \big(\mathrm{resp.,}~ C_c(G_{F,S},-):\mathrm{Ch}^+ (\Di_{\La, G_{F,S}})\lra
\mathrm{Ch}(Mod_{\La}) \big) \ea\] preserves homotopy, exact
sequences and quasi-isomorphisms, hence induces the following exact
derived functors
\[ \ba{c} \R\Ga_c(G_{F,S},-):\DD^b (\C_{\La, G_{F,S}})\lra
\DD(Mod_{\La}) \\
 \big(\mathrm{resp.,}~\R\Ga_c(G_{F,S},-):\DD^+ (\Di_{\La, G_{F,S}})\lra
\DD(Mod_{\La}) \big)\ea\] such that for $M$ in $\DD^b (\C_{\La,
G_{F,S}})$ or $\DD^+ (\Di_{\La, G_{F,S}})$, we have the following
exact triangle
\[\R\Ga_{c}(G_{F,S},M)\lra \R\Ga(G_{F,S},M)\lra \bigoplus_{v\in
S}\R\Ga(G_{v},M)\] in $\DD(Mod_{\La})$ and the following long exact
sequence
\[ \ba{c} \cdots \lra H^i_c(G_{F,S},M)\lra H^i(G_{F,S},M)
 \lra \displaystyle\bigoplus_{v\in
S}H^{i}(G_{v},M)\lra H^{i+1}_c(G_{F,S},M)\lra \cdots . \ea
\]
\ep

\bpf This is immediate from the definition of the cone. \epf

By \cite[Thm.\ 7.1.8(iii), Thm.\ 8.3.19]{NSW}, Proposition
\ref{inverse limit cochain} can be applied to $G_{F,S}$ and $G_v$,
where $v\in S_f$. For $v\in S_{\mathds{R}}$, $G_v$ is a finite group
of order 2, and so the finiteness hypothesis in Proposition
\ref{inverse limit cochain} is satisfied, so the conclusion also
holds in this case. The following analogous statement to Proposition
\ref{inverse limit cochain} for cohomology groups with compact
support will now follow from the definition of the cone and the long
exact sequence of cohomology groups in the preceding proposition.

\bp \label{cts compact cochain limit} The functor $C_c(G_{F,S},-)$
preserves direct limits in $\Di_{\La,G_{F,S}}$. Moreover, if $M =
\plim_{\al} M_{\al}$ is an object in $\C_{\La,G_{F,S}}$, where each
$M_{\al}$ is finite, then we have the following isomorphism
\[ C_{c}(G_{F,S},M) \cong \plim_{\al} C_{c}(G_{F,S},M_{\al}) \]
of complexes and isomorphisms
\[ H^{i}_c(G_{F,S},M) \cong \plim_{\al} H^i_c(G_{F,S},M_{\al})\]
of cohomology groups. \ep

\bl \label{global} We have
\begin{equation*}  H^{j}_c(G_{F,S},\Qp/\Zp(1))\cong
\begin{cases} \Qp/\Zp & \text{if $j=3$},\\  0 & \text{if $j>3$.}
\end{cases} \end{equation*}
\el

\bpf By the long exact sequence of Poitou-Tate \cite[8.6.13]{NSW},
we have the following exact sequence
\[ H^{2}(G_{F,S}, \Z/p^{n}\Z(1))\lra \bigoplus_{v\in S}H^{2}(G_{v},
\Z/p^{n}\Z(1)) \lra \Z/p^{n}\Z \lra 0 \] and an isomorphism \[
H^3(G_{F,S}, \Z/p^{n}\Z(1)) \stackrel{\mathrm{res}}{\lra}
\bigoplus_{v\in S_{\mathds{R}}}\hat{H}^{3}(G_{v}, \Z/p^{n}\Z(1)).\]
By the definition of continuous cochains with compact support and
the fact that $\cd_p(G_v) =2$ for $v\in S_f$, we have $
H^{3}_{c}(G_{F,S}, \Z/p^{n}\Z(1)) \cong \Z/p^{n}\Z$. The remainder
of the lemma will then follow from a similar argument to that in
Lemma \ref{local class field theory}. \epf

Let $M$ be a bounded complex in $\C_{\La,G_{F,S}}$. For each $v\in
S$, we define a morphism $\cup_v$ of complex of $\Zp$-modules to be
\[ C\big(G_v, M^{\vee}(1)\big)\ot_{\La} C\big(G_v, M) \lra C\big(G_v,
\Qp/\Zp(1)\big) \lra C_c\big(G_{F,S},\Qp/\Zp(1)\big)[1],
\] where the first map is the total cup product, and the second is
the natural morphism arising from the definition of the cone. By
Lemma \ref{global}, we have a quasi-isomorphism $\Qp/\Zp[-3]
\stackrel{i}{\lra}\tau_{\geq 3}C_c(G_{F,S}, \Qp/\Zp(1))$ of
complexes of $\Zp$-modules. Since $\Qp/\Zp$ is an injective
$\Zp$-module, the map $i$ has a homotopy inverse. We shall fix one
such map
\[r:\tau_{\geq 3}C_c(G_{F,S}, \Qp/\Zp(1))\lra \Qp/\Zp[-3], \] and this
induces the following morphism
\[  \vartheta : C_c(G_{F,S},\Qp/\Zp(1))\lra \tau_{\geq
3}C_c(G_{F,S}, \Qp/\Zp(1))\stackrel{r}{\lra} \Qp/\Zp[-3]\] of
complexes of $\Zp$-modules. Combining this with $\cup_v$, we obtain
a morphism
\[ C\big(G_v, M^{\vee}(1)\big)\ot_{\La} C\big(G_v, M) \stackrel{\cup_v}\lra C_c\big(G_{F,S},
\Qp/\Zp(1)\big)[1] \stackrel{\vartheta[1]}{\lra} \Qp/\Zp[-2]
\] of complexes of $\Zp$-modules. For $v\in S_f$,
this is essentially the morphism constructed in
Section 3.1, which will give the Tate local duality as in Theorem
\ref{Tate local adic}. We also have the following.

\bt \label{Tate local real 2} Let $p=2$, and let $v\in
S_{\mathds{R}}$. For a bounded complex $M$ of objects in
$\C_{\La,G_v}$, we have the following isomorphism
\[ \widehat{\R\Ga}(G_{v},M)\lra
\R\Hom_{\Zp}\Big(
\widehat{\R\Ga}\big(G_{v},M^{\vee}(1)\big),\Q_2/\Z_2\Big)[-2]
\] in $\DD(Mod_{\La})$. \et

\bpf By a similar argument to that in Theorem \ref{Tate local adic},
it suffices to consider a finite module $M$, and the conclusion then
follows from \cite[Thm.\ 7.2.17]{NSW}. \epf

For ease of notation, we shall write $P(G_{F,S}, -)$ for
$\bigoplus_{v\in S}C(G_v, -)$. We now define a morphism $\cup_S$ by
 \[
P\big(G_{F,S}, M^{\vee}(1)\big)\ot_{\La} P(G_{F,S}, M)
\stackrel{\nu}{\lra} \bigoplus_{v\in S}C_c\big(G_{F,S},
\Qp/\Zp(1)\big) \stackrel{\sum}{\lra} C_c\big(G_{F,S},
\Qp/\Zp(1)\big),\] where $\nu (a_S, b_S) = (a_v\cup_v b_v)_{v\in
S}$.

We now construct the total cup products for the compactly supported
cochain groups. Since these are defined as cones, it follows from
Lemma \ref{cone tech} that there are two morphisms
\[ \cup_0, \cup_1 : C_{c}(G_{F,S},M^{\vee}(1))\ot_{\La}C_c(G_{F,S},M)\lra
C_{c}(G_{F,S},\Qp/\Zp(1)) \] of complexes of $\Zp$-modules given by
\[ \ba{l} (a,a_S)\cup_{0}(b,b_{S}) = (a\cup b,
        (-1)^{\bar{a}}\mathrm{res}_{S}(a)\cup_S b_{S}) \\
    (a,a_{S}) \cup_1 (b,b_S) = (a\cup b, a_{S}\cup_S \mathrm{res}_{S}(b))    \ea \]
        where $\cup$ is the total cup product
\[C(G_{F,S},M^{\vee}(1))\ot_{\La}C(G_{F,S},M)\lra
C(G_{F,S},\Qp/\Zp(1)).\]  The morphisms $\cup_0$ and $\cup_1$ induce
the following morphisms

\[ \ba{c}  \cup_{c} : C(G_{F,S},M^{\vee}(1))\ot_{\La}C_{c}(G_{F,S},M)\lra
C_{c}(G_{F,S},\Qp/\Zp(1)) \\
_{c}\cup : C_{c}(G_{F,S},M^{\vee}(1))\ot_{\La}C(G_{F,S},M)\lra
C_{c}(G_{F,S},\Qp/\Zp(1)) \ea\] of complexes of abelian groups which
are given by the following respective formulas (see also
\cite[5.3.3.2, 5.3.3.3]{Ne})
\[ \ba{l} a\cup_{c}(b,b_{S}) = (a\cup b,
        (-1)^{\bar{a}}\mathrm{res}_{S}(a)\cup_S b_{S})\\
        (a,a_{S})\ _{c}\cup b = (a\cup b, a_{S}\cup_S \mathrm{res}_{S}(b)).\ea \]
All of these fit into the following diagram
\[ \SelectTips{eu}{} \xymatrix @C=0.7in{
  C_c(G_{F,S}, M^{\vee}(1))\ot_{\La}C_c(G_{F,S}, M) \ar[d]_{} \ar[r]^{} &
  C(G_{F,S}, M^{\vee}(1))\ot_{\La} C_c(G_{F,S}, M) \ar[d]^{\cup_{c}} \\
  C_c(G_{F,S}, M^{\vee}(1))\ot_{\La}C(G_{F,S}, M) \ar[r]^-{_c\cup} & C(G_{F,S}, \Qp/\Zp(1))   } \]
which is commutative up to homotopy by Lemma \ref{cone tech}. Also,
the following diagrams

\[  \SelectTips{eu}{} \entrymodifiers={!! <0pt, .8ex>+} \xymatrix @C=0.3in{
  C(G_{F,S}, M^{\vee}(1))\ot_{\La}\big(P(G_{v}, M)[-1]\big) \ar[dd]_{} \ar[r]^{} &
  P\big(G_{F,S}, M^{\vee}(1)\big) \ot_{\La}
  \big(P(G_{F,S}, M)[-1]\big) \ar[d]^-{t} \\
   & \Big(P\big(G_{F,S}, M^{\vee}(1)\big)\ot_{\La}P(G_{F,S}, M)
   \Big) [-1]
   \ar[d]^-{\cup_{S}[-1]} \\
  C(G_{F,S}, M^{\vee}(1))\ot_{\La}C_c(G_{F,S}, M) \ar[r]^{\cup_c} & C_c\big(G_{F,S}, \Qp/\Zp(1)\big)   } \]

  \[  \SelectTips{eu}{} \entrymodifiers={!! <0pt, .8ex>+} \xymatrix @C=0.3in{
  P(G_{F,S}, M^{\vee}(1))[-1]\ot_{\La} C(G_{F,S}, M) \ar[dd]_{} \ar[r]^{}
  & P(G_{F,S}, M^{\vee}(1))[-1]\ot_{\La} P(G_{F,S}, M)
   \ar[d]^-{t'} \\
   &      \big(P(G_{F,S}, M^{\vee}(1))\ot_{\La} P(G_{F,S}, M)\big) [-1]  \ar[d]^-{\cup_S[-1]} \\
  C_c(G_{F,S}, M^{\vee}(1)) \ot_{\La} C(G_{F,S}, M)
      \ar[r]^-{_c\cup} & C_c\big(G_{F,S}, \Qp/\Zp(1)\big)   } \]
are commutative, where $t$ and $t'$ are the morphisms defined as in
Lemma \ref{translation lemma}. These in turn induce the following
morphism of exact triangles in $\mathbf{K}(Mod_{\La})$.

\[  \SelectTips{eu}{} \entrymodifiers={!! <0pt, .8ex>+} \xymatrix @C=0.7in{
  P(G_{F,S}, M)[-1] \ar[d]_{} \ar[r]^{} &
  \displaystyle\bigoplus_{v\in S}\Hom_{\Zp}\Big(C(G_{v}, M^{\vee}(1)),C_c\big(G_{F,S}, \Qp/\Zp(1)\big)\Big)[-3] \ar[d]  \\
  C_c(G_{F,S}, M) \ar[d]_{} \ar[r]^{} &
  \Hom_{\Zp}\Big(C(G_{F,S}, M^{\vee}(1)), C_c\big(G_{F,S}, \Qp/\Zp(1)\big)\Big)[-3] \ar[d]  \\
  C(G_{F,S}, M) \ar[r]^{} & \Hom_{\Zp}\Big(C_c(G_{F,S}, M^{\vee}(1)), C_c\big(G_{F,S}, \Qp/\Zp(1)\big)\Big)[-3]   }
  \]

\noindent Combining this with the morphism \[ \vartheta:
C_c(G_{F,S},\Qp/\Zp(1))\lra \tau_{\geq 3}C_c(G_{F,S},
\Qp/\Zp(1))\stackrel{r}{\lra} \Qp/\Zp[-3],\] we obtain the following
morphism of exact triangles
\[  \SelectTips{eu}{} \entrymodifiers={!! <0pt, .8ex>+} \xymatrix @C=0.7in{
  \displaystyle\bigoplus_{v\in S}\R\Ga(G_{v}, M)[-1] \ar[d]_{} \ar[r]^{} &
  \displaystyle\bigoplus_{v\in S}\R\Hom_{\Zp}\big(\R\Ga(G_{v}, M^{\vee}(1)), \Qp/\Zp\big)[-3] \ar[d]  \\
  \R\Ga_c(G_{F,S}, M) \ar[d]_{} \ar[r]^{} &
  \R\Hom_{\Zp}\big(\R\Ga(G_{F,S}, M^{\vee}(1)), \Qp/\Zp\big)[-3] \ar[d]  \\
  \R\Ga(G_{F,S}, M) \ar[r]^{} & \R\Hom_{\Zp}\big(\R\Ga_c(G_{F,S}, M^{\vee}(1)), \Qp/\Zp\big)[-3]   }
  \] in $\DD(Mod_{\La})$.

\bt \label{Poitou-Tate adic} For any bounded complex $M$ in
$\C_{\La,G_{F,S}}$, the above morphism of exact triangles is an
isomorphism. \et

{\pf} : The top morphism is an isomorphism by Theorem \ref{Tate
local adic}. It remains to show that the middle morphism is an
isomorphism. By a similar argument (using Proposition \ref{cts
compact cochain limit} for the limiting argument for the compactly
supported cohomology) to that of Theorem \ref{Tate local adic}, we
can reduce to the case that $M$ is a single finite module. The
conclusion then follows from the usual Poitou-Tate duality (cf.
\cite[8.6.13]{NSW}). $\Box$

\br Theorem \ref{Tate local adic} and Theorem \ref{Poitou-Tate adic}
are stated in \cite{FK} for the case that $\La$ is a profinite ring
with a basis of neighborhoods consisting of powers of the Jacobson
radical of $\La$, and $M$ is a finitely generated projective
$\La$-module. \er

\section{Iwasawa modules} \label{Iwasawa modules}

In this section, we will introduce certain modules over an Iwasawa
algebra. The next two paragraphs will introduce some notations which
will be adhered to throughout this section.

Fix a prime $p$. Let $R$ be a commutative pro-$p$ ring with a
directed fundamental system $\mathcal{I}$ of neighborhoods of zero
consisting of open ideals. Let $G$ and $\Ga$ be two profinite groups
such that there is continuous homomorphism $\pi : G\lra \Ga$ of
profinite groups. Set $\La = R\llbracket \Ga\rrbracket$. We now
describe the natural profinite topology on $\La$ (see \cite[Sect.\
5.3]{RZ}). Let $\U$ be the collection of open normal subgroups of
$\Ga$, and consider the following family of two-sided ideals:
\[\ai\La + I(U), ~\ai\in\mathcal{I}, ~U\in\U.\] Here $I(U)$ denotes the
kernel of the map $\La \tha R[\Ga/U]$. We take these ideals as a
fundamental system of neighborhoods of zero.

We have a map $\iota : \La \ra \La$ which sends $\ga$ to $\ga^{-1}$.
Note that this is only a homomorphism of $R$-modules. It is a ring
homomorphism if and only if $\Ga$ is abelian. Denote by \[\rho =
\rho_{\Ga} : G \stackrel{\pi}{\lra} \Ga \sbs\La^{\times}\] the
tautological one-dimensional representation of $G$ over $\La$.

\br In most situations, the ring $R$ is usually a commutative
complete Noetherian local ring with finite residue field of
characteristic $p$, and the group $\Ga$ is a compact $p$-adic Lie
group. However, despite motivated by the above situation, we shall
consider the theory in more generality. \er

\subsection{Induced modules}

For a given $\La$-module $M$, we define a $\La^{\circ}$-module
$M^{\iota}$ by the formula $ m \cdot_{\iota} \la : = \iota(\la)m$
for $\la\in\La, m\in M$. Similarly, if $N$ is a
$\La^{\circ}$-module, we define a $\La$-module, which is also
denoted as $N^{\iota}$, by $\la \cdot_{\iota} m : = m\iota(\la)$.

We shall prove the following lemma. Let $A$ and $B$ be two rings,
and suppose that $M$ has a left $A$-action and right $B$-action. We
say that the actions of $A$ and $B$ are balanced if for every $a\in
A,b\in B$ and $x\in M$, we have $a(xb) = (ax)b$.

\bl  $(a)$ If $M$ is a $\La[G]$-module, then $M^{\iota}$ is a
$\La^{\circ}[G]$-module.

$(b)$ If $M$ is a $\La[G]$-$\La$-module $($not necessarily
balanced$)$, then $M^{\iota}$ is a $\La^{\circ}[G]$-$\La$-module
$($not necessarily balanced$)$.
 \el

\bpf (a) Let $g\in G, \la\in\La$ and $m\in M^{\iota}$. Then we have
\[ (gm)\cdot_{\iota}\la = \iota(\la)gm = g(\iota(\la)m) =
g(m\cdot_{\iota}\la).\]

(b) Similar argument as above. \epf

\medskip
For a given $U\in\U$ and a given $R[G]$-module $M$, we define two
$\La[G]$-$\La$-modules as follows:
\[\ba{c}  _{U}M = \Hom_{R}(R[\Ga/U],M) \\
     M_{U} = R[\Ga/U]^{\iota}\ot_{R}M, \\ \ea \]
where $G$ acts on $R[\Ga/U]$ via $\rho_{\Ga/U}$ and $\La$ acts on
$R[\Ga/U]$ via the canonical projection $\La\tha R[\Ga/U]$. Note
that the $\La[G]$-$\La$-modules defined above are balanced as
$\La$-$\La$-modules. They are balanced as $\La[G]$-$\La$-modules if
$\Ga /U$ is abelian.

Let $V\in\U$ with $U\sbs V$. Then there is a canonical surjection
$\mathrm{pr}:R[\Ga/U]\tha R[\Ga/V]$ and a map
$\mathrm{Tr}:R[\Ga/V]\lra R[\Ga/U]$ given by
\[ gU \mapsto\sum_{v\in V/U}gvU.\] These in turn induce the
following maps.
\[ \ba{c} \mathrm{pr}^{*}: ~_{V}M\lra\, _{U}M \\
          \mathrm{pr}_{*}: M_{U}\lra M_{V}\\
          \mathrm{Tr}^{*}: ~_{U}M\lra\, _{V}M \\
          \mathrm{Tr}_{*}: M_{V}\lra M_{U}
          \ea \]
Denote by $\de_{\be} : G/U \ra \Z$ the Kronecker delta-function
\begin{equation*}  \de_{\be}(\be')\cong
\begin{cases} 1 & \text{if $\be=\be'$},\\  0 & \text{if $\be\neq \be'$.}
\end{cases} \end{equation*}
The next two lemmas then follow from a straightforward calculation.

\bl We have the following isomorphism of $R[G]$-modules
   \[ \ba{c}   M_{U}\stackrel{\sim}{\lra}\, _{U}M \\
           \displaystyle \sum_{\be\in G/U}\be\ot x_{\be}\mapsto \sum_{\be\in
             G/U}x_{\be}\de_{\be} \ea \] which is functorial in $M$.
Moreover, if $V$ is another open normal subgroup of $G$ such that
$U\sbs V$, then the isomorphism fits into the following commutative
diagrams.
\[\xymatrix{
  M_{U} \ar[d]_{\mathrm{pr}_{*}} \ar[r]^{\sim} & {}_U M \ar[d]^{\mathrm{Tr}^{*}}
  && M_{V} \ar[d]_{\mathrm{Tr}_{*}} \ar[r]^{\sim} & _{V}M \ar[d]^{\mathrm{pr}^{*}}\\
  M_{V} \ar[r]^{\sim} & _{V}M  &&M_{U} \ar[r]^{\sim} & _{U}M }\]
   \el

\bl We have the following equalities of $\La^{\circ}[G]$-modules$:$
\[ \ba{c}
 (_{U}M)^{\iota} = \Hom_{R}(R[\Ga/U]^{\iota},M),\\
 (M_{U})^{\iota} = R[\Ga/U]\ot_{R}M. \ea \] \el

\medskip Let $M$ be an $R[G]$-module. We define two
$\La[G]$-$\La$-modules as follows:
\[ \ba{c}
F_{\Ga}(M) = \ilim_{U\in\U} {}_{U}M,\\ \F_{\Ga}(M) = \plim_{U\in\U}
M_{U},  \ea \] where the transition maps are induced by the
surjections $R[\Ga/U]\tha R[\Ga/V]$ for $U\sbs V$. Note that the
$\La[G]$-$\La$-modules defined above are balanced as
$\La$-$\La$-modules. They are balanced as $\La[G]$-$\La$-modules if
and only if $\Ga$ is abelian. One easily sees from Lemma 3.1.3 that
\[ \ba{c}F_{\Ga}(M)^{\iota} = \ilim_{U\in\U} \Hom_{R}(R[\Ga/U]^{\iota},
M)~\mathrm{and}\\  \F_{\Ga}(M)^{\iota} = \plim_{U\in\U}
(R[\Ga/U]\ot_{R} M). \ea
\]
We also have the following description of $F_{\Ga}(A)$, when $A$ is
an object of $\Di_{R,G}$.

\bl \label{F(M)2} If $A$ is an object of $\Di_{R,G}$, then
$F_{\Ga}(A)$ is an object of $\Di_{\La,G}$ and
\[F_{\Ga}(A) \cong \Hom_{R,\cts}(\La,A).\]

\noindent Similarly, we have
\[F_{\Ga}(A)^{\iota} \cong
\Hom_{R,\cts}(\La^{\iota},A).\]

If $\{A_{\al}\}$ is a direct system of objects in $\Di_{R,G}$, then
we have isomorphisms
\[ F_{\Ga}(A) \cong \ilim_{\al} F_{\Ga}(A_{\al}) ~\big(resp.,\
F_{\Ga}(A)^{\iota} \cong \ilim_{\al} F_{\Ga}(A_{\al})^{\iota}\big)\]
in $\Di_{\La,G}$ $($resp., in $\Di_{\La^{\circ},G})$. \el

\bpf By Lemma \ref{abstract is cts}(3), for each $U\in\U$, we have
\[ \Hom_{R}(R[\Ga/U],A) = \Hom_{R,\cts}(R[\Ga/U],A). \]
Therefore, the lemma will now follow from \cite[Prop.\ 5.1.4]{RZ}
\epf

We would like to have a description of $\F_{\Ga}(T)$, when $T$ is an
object in $\C_{R,G}$. Before we can do this, we shall recall the
notion of a complete tensor product from \cite{RZ}. Let $M$ be an
object in $\C_{\La,G}$, and let $N$ be an object in $\C_{R,G}$. The
completed tensor product of $M$ and $N$ is taken to be
\[  M\otw_{R}N = \plim_{U,V}M/U\ot_{R}N/V,\]
where $U$ (resp.,\ $V$) runs through the open $\La[G]$-submodules of
$M$ (resp.,\ open $R[G]$-submodules of $N$).

\bl \label{completed tensor product2} Let $M$ be an object in
$\C_{\La,G}$ and $N$ be an object in $\C_{R,G}$.  Then the completed
tensor product $M \otw_{R} N$ is an object of $\C_{\La,G}$, and
coincides with the usual tensor product if $N$ is a finitely
generated $R$-module. Moreover, as a functor, the completed tensor
product is right exact $($in both variables$)$ and preserves inverse
limits. \el

\bpf It follows from \cite[Lemma 7.7.2]{Wil} that $M\otw_{R}N$ is a
compact $\La$-module. By a similar argument to that used in the
proof of that lemma, we have that the $G$-action is continuous. \epf

We are now in position to describe $\F_{\Ga}(T)$.

\bl \label{F(M)} If $T$ is an object of $\C_{R,G}$, then
$\F_{\Ga}(T)$ is isomorphic to $\La^{\iota}\otw_{R}T$ and
$\F_{\Ga}(T)^{\iota}$ is isomorphic to $\La\otw_{R}T$. If
$\{T_{\al}\}$ is an inverse system of objects in $\C_{R,G}$ such
that $T \cong \plim_{\al} T_{\al}$, then we have isomorphisms
\[ \F_{\Ga}(T) \cong \plim_{\al} \F_{\Ga}(T_{\al})~\big(resp.,\
\F_{\Ga}(T)^{\iota} \cong \plim_{\al} \F_{\Ga}(T_{\al})^{\iota}\big)
\] in $\C_{\La,G}$ $($resp.,\ in $\C_{\La^{\circ},G})$. \el

\bpf We have
\[ \F_{\Ga}(T) = \plim_{U}(R[\Ga/U]^{\iota}\ot_R T) =
\plim_{U} \big(R[\Ga/U]^{\iota}\otw_R T\big) \cong
\big(\plim_{U}R[\Ga/U]^{\iota}\big)\otw_R T \cong \La^{\iota}\otw_R
T. \] Suppose $T \cong \plim_{\al} T_{\al}$ in $\C_{R,G}$. Then
\[ \F_{\Ga}(T) \cong \La^{\iota}\otw_R
T \cong \plim_{\al} \La^{\iota}\otw_R T_{\al} \cong \plim_{\al}
\F_{\Ga}(T_{\al}). \] \epf

\medskip
As a conclusion to the subsection, we record the following duality
relation between the modules we have defined.

\bp \label{DFM} Let $T$ be an object in $\C_{R,G}$. Then we have
isomorphisms
\[\F_{\Ga}(T)^{\vee} \cong
F_{\Ga}(T^{\vee})^{\iota}~\Big(resp.,\ (\F_{\Ga}(T)^{\iota})^{\vee}
\cong F_{\Ga}(T^{\vee})\Big)\] in $\Di_{\La^{\circ},G}$ $($resp., in
$\Di_{\La,G})$.  \ep

\bpf We will prove the first isomorphism, the second will follow
from a similar argument. This follows by the following calculations:
 \[\begin{array}{rl}
          \F_{\Ga}(T)^{\vee} \cong & \Hom_{\Zp}(\La^{\iota}\otw_R T, \Qp/\Zp) ~(\mbox{by Lemma \ref{F(M)}}) \\
          \cong & \Hom_{R,\cts}\big(\La^{\iota},
          \Hom_{\Zp,\cts}(T, \Qp/\Zp)\big) ~(\mbox{by \cite[Prop.\ 5.5.4(c)]{RZ}})\\
          \cong &  F_{\Ga}(T^{\vee})^{\iota} ~(\mbox{by Lemma \ref{F(M)2}}).\ea
          \] \epf

\subsection{Shapiro's lemma} \label{Shapiro section}

As before, $R$ denotes a commutative pro-$p$ ring. Let $G$ be a
profinite group. Fix a closed normal subgroup $H$ of $G$ and write
$\Ga = G/H$. Let $\pi :G\lra \Ga$  be the canonical quotient map. We
identify $\U$ as the collection of open normal subgroups of $G$
containing $H$. Therefore, in this context, for each $U\in\U$, and
an $R[G]$-module $M$, we have
\[\ba{c}  _{U}M = \Hom_{R}(R[G/U],M), \\
     M_{U} = R[G/U]^{\iota}\ot_{R}M. \\ \ea \]
We will apply Shapiro's lemma to see that the direct limits and
inverse limits of cohomology groups over every intermediate field
$F_{\al}$ can be viewed as cohomology groups of certain
$\Lambda$-modules. The results in this section can be found in
\cite[8.2.2, 8.3.3-5, 8.4.4.2]{Ne}.

\bl \label{Shapiro discrete} Let $U$ be an open normal subgroup of
$G$, and let $N$ be a bounded below complex of objects of
$\Di_{R,G}$. Then we have a quasi-isomorphism
 \[ C(G,{}_{U}N) \stackrel{\sim}{\lra} C(U, N)
\] of complexes of $\La$-modules. \el

\bpf We first prove the lemma in the case that $N$ is an object of
$\Di_{R,G}$. Then we may write $N = \ilim_{\al}N_{\al}$, where
$N_{\al}$ is a finite object of $\Di_{R,G}$. The usual Shapiro's
lemma holds for such modules. Also, we note that $_UN \cong
\ilim_{\al}\,_U(N_{\al})$. Hence, we have \[ C(G, \,_UN) = C\Big(G,
\ilim_{\al}\,_U(N_{\al})\Big) \cong \ilim_{\al}C\big(G,
\,_U(N_{\al})\big) \stackrel{\mathrm{sh}}{\lra} \ilim_{\al}C(U,
N_{\al}) = C(U, N) \] which gives the required conclusion for the
case that $N$ is an object of $\Di_{R,G}$. For the case that $N$ is
a bounded below complex of objects of $\Di_{R,G}$, one can prove
this by a spectral sequence argument as used in Lemma \ref{bounded
below qis}. \epf

Recall that if $A$ is a complex in $\Di_{R,G}$, then $F_{\Ga}(A) =
\ilim_{U\in\U} \,_U A$ is a complex in $\Di_{\La,G}$ by Lemma
\ref{F(M)2}. We then have the following proposition.

\bp \label{Shapiro1} Let $A$ be a bounded below complex of objects
of $\Di_{R,G}$. Then the composite morphism
\[ C(G, F_{\Ga}(A)) \stackrel{\sim}{\lra} \ilim_{U\in\U}C(G,
\, _{U}A) \stackrel{\mathrm{sh}}{\lra} \ilim_{U\in\U}C(U, A)
\stackrel{\mathrm{res}}{\lra} C(H,A)\] is a quasi-isomorphism of
complexes of $\La$-modules. In other words, we have an isomorphism
\[ \R\Ga\big(G, F_{\Ga}(A)\big) \stackrel{\sim}{\lra} \R\Ga(H,A) \] in
$\DD(Mod_{\La})$. \ep

The next result will give a Shapiro-type relation for cohomology
groups of objects (and complexes of objects) in $\C_{R,G}$.

\bl \label{Shapiro compact} Let $U$ be an open normal subgroup of
$G$. Then for any bounded complex $M$ in $\C_{R,G}$, we have a
quasi-isomorphism
 \[ C(G, M_U) \stackrel{\sim}{\lra} C(U, M)
\] of complexes of $\La$-modules. \el

\bpf By the same argument as that in Lemma \ref{Shapiro discrete},
it suffices to consider the case when $M$ is an object of
$\C_{R,G}$. Then we have $M = \plim_{\al}M_{\al}$, where $M_{\al}$
is a finite object in $\C_{R,G}$. Note that $M_U\cong
\plim_{\al}(M_{\al})_U$. Then we have morphisms \[ C(G, M_U) \cong
\plim_{\al}C\big(G, (M_{\al})_U\big) \stackrel{\mathrm{sh}}{\lra}
\plim_{\al}C(U, M_{\al}) \cong C(U, M)
\] which induce a morphism
\[\xymatrix{
  {\plim_{\al}}^{(i)} H^j\big(G,(M_{\al})_U\big) \Longrightarrow H^{i+j}(G, M_U)  \ar[d]^{} \\
  {\plim_{\al}}^{(i)} H^j\big(U,M_{\al}\big) \Longrightarrow H^{i+j}(U, M)   }
  \] of convergent spectral sequences.
Since $M_{\al}$ is finite, the usual Shapiro's lemma implies that
\[H^j\big(G,(M_{\al})_U\big)\cong H^j\big(U,M_{\al}\big)\] is an isomorphism. This
in turn implies that
\[ {\plim_{\al}}^{(i)} H^j\big(G,(M_{\al})_U\big) \cong {\plim_{\al}}^{(i)} H^j\big(U,M_{\al}\big).\]
By the convergence of the spectral sequences, we have isomorphisms
\[H^{n}(G, M_U) \cong  H^n(U,M), \] as required. \epf

Since inverse limits are not necessarily exact, we cannot always
view inverse limits of cohomology groups over every intermediate
field $F_{\al}$ as cohomology groups of certain $\Lambda$-modules in
general. However, we can say something if we impose an extra
assumption on $G$.

\bp \label{Shapiro2} Let $M$ be a bounded complex of objects in
$\C_{R,G}$. Then we have the following isomorphism
\[ C(G,\F_{\Ga}(M)) \stackrel{\sim}{\lra} \plim_{U\in\U}C(G, M_{U}) \] of complexes of $\La$-modules.
Furthermore, if $H^{m}(G,N)$ is finite for all finite discrete
$\La$-modules $N$ with a $\La$-linear continuous $G$-action and all
$m\geq 0$, then we have
\[H^{j}(G,\F_{\Ga}(M)) \cong \plim_{U\in\U}H^{j}(U,M).\]\ep

\bpf As before, it suffices to consider the case when $M$ is an
object of $\C_{R,G}$. Write $M = \plim_{\al} M_{\al}$, where each
$M_{\al}$ is a finite object in $\C_{R,G}$. By Lemma \ref{F(M)}, we
have a continuous isomorphism
\[ \F_{\Ga}(M)\cong \plim_{\al} \F_{\Ga}(M_{\al}) \cong  \plim_{\al,U}
(M_{\al})_{U} .\]  The second assertion now follows from Proposition
\ref{inverse limit cochain} and Lemma \ref{Shapiro compact}. \epf

\subsection{Iwasawa setting}

We now apply the discussion in Subsection \ref{Shapiro section} to
the arithmetic situation. Let $F_{\infty}$ be a Galois extension of
$F$ which is contained in $F_{S}$. Write $H =
\mathrm{Gal}(F_{S}/F_{\infty})$, and write $\Ga =
\mathrm{Gal}(F_{\infty}/F)$. Let $\U$ denote the collection of open
normal subgroups of $G_{F,S}$ containing $H$. For each $U\in\U$, we
let $F_{U} = (F_{S})^{U}$, and define $S_U$ to be the set of primes
in $F_U$ above $S$.  As before, we write $\La =
R\llbracket\Ga\rrbracket$, where $R$ is a commutative pro-$p$ ring.
The following lemma is immediate from the discussion in the
preceding subsection.

\bl \label{globalshapiro}  Let $T$ be a bounded complex of objects
in $\C_{R,G_{F,S}}$. Then we have the following isomorphisms
\[ \ba{c}
H^{j}(G_{F,S},\F_{\Ga}(T)) \cong \plim_{U}H^j(G_{F,S}, T_{U}) \cong \plim_{U}H^{j}(G_{F_U,S_U},T), \\
H^{j}(G_{F,S},F_{\Ga}(T^{\vee})) \cong
\ilim_{U}H^j(G_{F,S},\,_{U}T^{\vee})
 \cong\ilim_{U}H^{j}(G_{F_U,S_U},T^{\vee}) \cong H^{j}(\mathrm{Gal}(F_S/F_{\infty}),T^{\vee}).
 \ea\]
\el

Let $v\in S_f$. Fix an embedding $F^{\mathrm{sep}}\hra
F_v^{\mathrm{sep}}$, which induces a continuous group monomorphism
\[ \al = \al_v : G_v \hra G_F,\] where $G_F =
\mathrm{Gal}(F^{\mathrm{sep}}/F)$.  Let $X$ be a finite discrete
$R[G_{F}]$-module. For a finite Galois extension $F'$ of $F$, write
$U=\mathrm{Gal}(F^{\mathrm{sep}}/F')$ and $X_{U} =
R[G_{F}/U]\ot_{R}X$. The embedding $F^{\mathrm{sep}}\hra
F_v^{\mathrm{sep}}$ determines a prime $v'$ of $F'$ above $v$ such
that $F'_{v'}$ is a finite Galois extension of $F_v$ and $G_{v'}
:=\mathrm{Gal}(F^{\mathrm{sep}}_v/F'_{v'}) = \al^{-1}(U)$.

Fix coset representatives $\s_i\in G_{F}$ of
\[ G_{F}/U = \bigcup_i\s_i \al(G_v/G_{v'}). \] Then the set of
distinct primes in $F'$ above $v$ is given by the (finite)
collection $\{\s_{i}(v')\}$. Then by \cite[8.1.7.6, 8.5.3.1]{Ne}, we
have a quasi-isomorphism
\[ C(G_v, X_U)\stackrel{\sim}{\lra} \bigoplus_i C(G_{\s_{i}(v')}, X)   \]  and isomorphisms
\[ H^n(G_v , X_U) \cong \bigoplus_i H^n(G_{\s_{i}(v')}, X)\] of cohomology
groups for $n \geq 0$.

Now suppose that $p=2$ and $F$ is a number field with at least one
real prime. Let $v\in S_{\mathds{R}}$.  Let $F'$ be a finite Galois
extension of $F$, and retain the above notations. Then the primes in
$F'$ above $v$ are either all real or all complex. We first consider
the case when all the primes above $F'$ are all real. Then $G_v$
acts trivially on $X_{U}$, and we have
\[ C(G_v, X_U) =  C(G_v, X)^{|G_F/U|}\] which is precisely $ \bigoplus_{\tau} C(G_{\tau(v')},
X)$, where $v'$ is a prime above $v$ and $\tau$ runs through a set
of coset representatives for $G_F/U$. If all the primes above $F'$
are complex, it remains to show that $\hat{H}^i(G_v, X_U)=0$ for all
$i$. Since $G_v$ is cyclic (of order 2), we are reduced to showing
this for $i=1,2$, which follows from Shapiro's lemma in the usual
sense (since these are usual cohomology groups)

 We shall apply the above discussion to finite discrete
$R[G_{F,S}]$-modules, which we view as $R[G_F]$-modules via the
canonical quotient map $G_F\tha G_{F,S}$. By the compatibility of
limits and the groups of continuous cochains, we can apply the above
results to objects in $\Di_{R,G_{F,S}}$ or $\C_{R,G_{F,S}}$.

\bl \label{localshapiro2}  Let $T$ be a bounded complex of objects
in $\C_{R,G_{F,S}}$. Then for $v\in S_f$, we have the following
isomorphisms
\[ \ba{c}
H^{j}(G_{v}, \F_{\Ga}(T)) \cong \plim_{U}H^j(G_{v}, T_{U}) \cong
\plim_{U}\bo_{w|v}H^{j}(G_{w}, T),\\
H^{j}(G_{v},F_{\Ga}(T^{\vee})) \cong
\ilim_{U}H^j(G_{v},\,_{U}T^{\vee}) \cong
\ilim_{U}\bo_{w|v}H^{j}(G_{w}, T^{\vee}). \ea \] The same conclusion
holds for the case when $p=2$ and $v\in S_{\mathds{R}}$, if we
replace the cohomology groups by the completed cohomology groups as
defined in Subsection \ref{Tate cohomology groups}. \el

\medskip
We would like to derive an analogue of Shapiro's lemma for compactly
supported cohomology. Let $F'$ be a finite Galois extension of $F$
which is contained in $F_{S}$. Denote the set of primes of $F'$
above $S$ by $S'$. Let $X$ be a discrete $R[G_{F,S}]$-module. We
write $U=\mathrm{Gal}(F_{S}/F')$ and $X_{U} = R[G_{F,S}/U]\ot_{R}X$.
By the discussion in the previous subsection and the above, we have
the following diagram

\[\SelectTips{eu}{} \xymatrix@C=0.45in{
  C(G_{F,S},X_{U}) \ar[d]^{\mathrm{sh}} \ar[r]^-{}
& \bo_{v\in S}C(G_{v},X_{U}) \ar[r]^-{\sim}& \bo_{v\in
S}\bo_{v'|v}C(G_{v},R[G_{v}/G_{v'}]\ot_{R}X)
 \ar[d]^{\mathrm{sh}} \\
  C(G_{F',S'},X) \ar[rr]^-{} &&\bo_{v'\in S'}C(G_{v'},X)
}\] which commutes up to homotopy. By a similar argument to that in
\cite[8.1.7.2.1, 8.5.3.2]{Ne}, this in turn induces a
quasi-isomorphism (functorial in $X$)
\[\mathrm{sh}_{c} : C_{c}(G_{F,S}, X_{U}) \lra C_{c}(G_{F',S'}, X)
\]
which fits into the following commutative (up to homotopy) diagram
with exact rows.
\[\SelectTips{eu}{} \xymatrix @C=0.25in{
  0 \ar[r]^{} & \bo_{v\in S}C(G_{v},X_{U})[-1] \ar[d]_{\wr} \ar[r]^{} & C_{c}(G_{F,S}, X_{U}) \ar @{=}[d]
  \ar[r]^{} & C(G_{F,S}, X_{U}) \ar @{=}[d] \ar[r]^{} & 0 \\
  0 \ar[r]^{} & \bo_{v\in S}\bo_{v'|v}C(G_{v},R[G_{v}/G_{v'}]\ot_{R}X)[-1] \ar[d]_{\mathrm{sh}[-1]} \ar[r]^{}
& C_{c}(G_{F,S}, X_{U}) \ar[d]_{\mathrm{sh}_{c}}
  \ar[r]^{} & C(G_{F,S}, X_{U}) \ar[d]_{\mathrm{sh}} \ar[r]^{} & 0 \\
  0 \ar[r]^{} & \bo_{v'\in S'}C(G_{v},X)[-1] \ar[r]^{} & C_{c}(G_{F',S'}, X)
  \ar[r]^{} & C(G_{F',S'}, X) \ar[r]^{} & 0  } \]

Suppose that $F''\sbs F_S$ is another finite Galois extension of $F$
containing $F'$, and write $S''$ for the set of primes of $F''$
above $S$ and $V = \mathrm{Gal}(F_{S}/F'')$. Again by similar
arguments to that in \cite[8.5.3.4]{Ne}, we have the following
morphisms
\[ \ba{c} \mathrm{res}_{c} : C_{c}(G_{F',S'},X) \lra C_{c}(G_{F'',S''},X) \\
     \mathrm{cor}_{c} : C_{c}(G_{F'',S''},X) \lra C_{c}(G_{F',S'},X), \ea \]
which are functorial in $X$ and fit in the following diagrams, which
are commutative up to homotopy:
\[\SelectTips{eu}{}\xymatrix{
  0  \ar[r] & \oplus_{v'\in S'}C(G_{v'},X)[-1] \ar[d]^{\mathrm{res}[-1]} \ar[r] & C_{c}(G_{F',S'},X)
  \ar[d]^{\mathrm{res}_{c}}
  \ar[r] & C(G_{F',S'},X) \ar[d]^{\mathrm{res}} \ar[r] & 0  \\
  0  \ar[r] & \oplus_{v''\in S''}C(G_{v''},X)[-1] \ar[r] & C_{c}(G_{F'',S''},X)
  \ar[r] & C(G_{F'',S''},X)  \ar[r] & 0   } \]

\[\SelectTips{eu}{}\xymatrix{
  0  \ar[r] & \oplus_{v''\in S''}C(G_{v''},X)[-1] \ar[d]^{\mathrm{cor}[-1]} \ar[r] & C_{c}(G_{F'',S''},X)
  \ar[d]^{\mathrm{cor}_{c}}
  \ar[r] & C(G_{F'',S''},X) \ar[d]^{\mathrm{cor}} \ar[r] & 0  \\
  0  \ar[r] & \oplus_{v'\in S'}C(G_{v'},X)[-1] \ar[r] & C_{c}(G_{F',S'},X)
  \ar[r] & C(G_{F',S'},X)  \ar[r] & 0   } \]

\[\SelectTips{eu}{}\xymatrix{
  C_{c}(G_{F,S},X_{U}) \ar[d]_{\mathrm{Tr}_{*}} \ar[r]^{\mathrm{sh}_{c}} & C_{c}(G_{F',S'},X) \ar[d]_{\mathrm{res}_{c}}
  & C_{c}(G_{F,S},X_{V}) \ar[d]_{\mathrm{pr}_{*}}
  \ar[r]^{\mathrm{sh}_{c}} & C_{c}(G_{F'',S''},X) \ar[d]^{\mathrm{cor}_{c}} \\
  C_{c}(G_{F,S},X_{V}) \ar[r]^{\mathrm{sh}_{c}} & C_{c}(G_{F'',S''},X)
   & C_{c}(G_{F,S},X_{U}) \ar[r]^{\mathrm{sh}_{c}} & C_{c}(G_{F',S'},X)
     }\]

Since all the morphisms constructed above are functorial, they can
be extended to complexes. Hence, we may conclude the following.

\bp \label{compactshapiro} $(a)$ For a bounded below complex $A$ of
objects in $\Di_{R,G_{F,S}}$, the canonical morphism of complexes \[
C_{c}(G_{F,S},F_{\Ga}(A))\stackrel{\sim}{\lra}\ilim_{U,\mathrm{Tr}}C_{c}(G_{F,S},
\,_{U}A)
\] is an isomorphism.

$(b)$  Let $T$ be a bounded complex of objects in $\C_{R,G_{F,S}}$.
Then the canonical morphism of complexes
\[
C_{c}(G_{F,S},\F_{\Ga}(T))\stackrel{\sim}{\lra}\plim_{U,
\mathrm{pr}}C_{c}(G_{F,S},T_{U})
\] is an isomorphism and induces isomorphisms
\[ H^{j}_c(G_{F,S},\F_{\Ga}(T)) \cong \plim_{U, \mathrm{pr}}H^j_{c}(G_{F,S},T_{U})
\cong \plim_{U, \mathrm{cor}_c}H^j_{c}(G_{F_U,S_U},T)\] of
cohomology groups for $j\geq 0$. \ep

\subsection{Duality over extensions of a global/local field}

We retain the notation introduced in the previous subsection. Let
$F_{\infty}$ be a Galois extension of $F$ which is contained in
$F_{S}$. Write $H = \mathrm{Gal}(F_{S}/F_{\infty})$, and write $\Ga
= \mathrm{Gal}(F_{\infty}/F)$. As before, we write $\La =
R\llbracket\Ga\rrbracket$, where $R$ is a commutative pro-$p$ ring.

 Applying Theorem \ref{Poitou-Tate adic} and
Proposition \ref{DFM}, we obtain the following theorem. In the
theorem, we abuse notation and use $\R\Ga(G_v,M)$ to denote
$\widehat{\R\Ga}(G_v,M)$ for $v\in S_{\mathds{R}}$.

\bt \label{p-adic Tate duality}  Then, for a bounded complex $T$ in
$\C_{R,G_{F,S}}$, we have the following isomorphism of exact
triangles
\[ \SelectTips{eu}{} \entrymodifiers={!! <0pt, .8ex>+} \xymatrix @C=0.65in{
  \displaystyle\bigoplus_{v\in S}\R\Ga(G_{v}, \F_{\Ga}(T))[-1] \ar[d]_{} \ar[r]^-{\sim} &
  \displaystyle\bigoplus_{v\in S}\R\Hom_{\Zp}\big(\R\Ga(G_{v}, F_{\Ga}(T^{\vee})^{\iota}(1)), \Qp/\Zp\big)[-3] \ar[d]  \\
  \R\Ga_c(G_{F,S},  \F_{\Ga}(T)) \ar[d]_{} \ar[r]^-{\sim} &
  \R\Hom_{\Zp}\big(\R\Ga(G_{F,S}, F_{\Ga}(T^{\vee})^{\iota}(1)), \Qp/\Zp\big)[-3] \ar[d]  \\
  \R\Ga(G_{F,S},  \F_{\Ga}(T)) \ar[r]^-{\sim} & \R\Hom_{\Zp}\big(\R\Ga_c(G_{F,S}, F_{\Ga}(T^{\vee})^{\iota}(1)), \Qp/\Zp\big)[-3]   }
  \] in $\DD(Mod_{\La})$.  \et

We end by saying something about the situation over nonarchimedean
local fields. Let $F$ be a nonarchimedean local field of
characteristic not equal to $p$. Let $F_{\infty}$ be a Galois
extension of $F$ with Galois group $\Ga$. Write $G_E =
\mathrm{Gal}(F^{\mathrm{sep}}/E)$ for every Galois extension $E/F$.
Recall that by \cite[Thm.\ 7.1.8(i)]{NSW}, we have
$\mathrm{cd}_p(G_F) = 2$.

Let $T$ be a bounded complex of objects in $\C_{R,G_F}$. By
Proposition \ref{Shapiro1} and Proposition \ref{Shapiro2}, we have
\[ \ba{c}
C\big(G_F, F_{\Ga}(T^{\vee})\big) \stackrel{\sim}{\lra} \ilim
C(G_{F_{\al}}, T^{\vee}) \\
H^i\big(G_F, F_{\Ga}(T^{\vee})\big) \cong \ilim
H^i(G_{F_{\al}}, T^{\vee}) \cong H^i(G_{F_{\infty}},T^{\vee})\\
C\big(G_F, \F_{\Ga}(T)\big) \stackrel{\sim}{\lra} \plim
C(G_{F_{\al}}, T) \\
H^i\big(G_F, \F_{\Ga}(T)\big) \cong \plim H^i(G_{F_{\al}}, T),  \ea
\] where $F_{\al}$ runs through all finite Galois extension of
$F_{\infty}/F$. Applying Theorem \ref{Tate local adic} and
Proposition \ref{DFM}, we obtain the following.

\bt Let $T$ be a bounded complex of objects in $\C_{R,G_{F}}$. Then
we have the following isomorphism
\[ \R\Ga(G_{F},\F_{\Ga}(T)) \stackrel{\sim}{\lra}
\R\Hom_{\Zp}\Big(\R\Ga\big(G_{F},F_{\Ga}(T^{\vee})^{\iota}(1)\big),\Qp/\Zp\Big)[-2]
\]
in $\DD(Mod_{\La})$.  \et


\begin{thebibliography}{15}

\bibitem{CFKSV} J.\ Coates, T.\ Fukaya, K.\ Kato, R.\ Sujatha and O.\ Venjakob, The
$GL_2$ main conjecture for elliptic curves without complex
multiplication, \textit{Publ. Math. Inst. Hautes \'Etudes Sci.},
\textbf{101} (2005), 163-208.

\bibitem{F} J.\ Flood, Pontryagin duality for topological modules,
 \textit{Proc. Amer. Math. Soc.} vol. \textbf{75} (1979), No. 2,
 329-333.

\bibitem{FK} T.\ Fukaya and K.\ Kato, A formulation of
conjectures on $p$-adic zeta functions in noncommutative Iwasawa
theory, \textit{Amer. Math. Soc. Transl. Ser. 2} vol.
$\mathbf{219}$, (2006), 1-85.

\bibitem{RD} R.\ Hartshorne, \textit{Residues and Duality}, Lecture
Notes in Math. \textbf{20}, Springer-Verlag, Berlin, (1966).


\bibitem{Je} C.\ U.\ Jensen, \textit{Les Foncteurs D\'eriv\'es de
$\plim$ et leurs Applications en Th\'eorie des Modules}, Lecture
Notes in Math. \textbf{254}, Springer-Verlag, Berlin, (1972).

\bibitem{Ka} M.\ Kakde, The main conjecture of Iwasawa theory for totally real
fields, arXiv:1008.0142v2 [math.NT].

\bibitem{Lim-thesis} M. F. Lim, Duality over $p$-adic Lie extensions of global fields, Ph.D. thesis,
McMaster University, (2010).

\bibitem{LS} M. F. Lim and  R. Sharifi, Nekov\'a\v{r} duality over $p$-adic Lie extensions of global
fields. Preprint.

\bibitem{Ne} J.\ Nekov\'a\v{r}, Selmer complexes,
$Ast\acute{e}risque$ \textbf{310} (2006).


\bibitem{NSW} J.\ Neukirch, A.\ Schmidt and K.\ Wingberg,
\textit{Cohomology of Number Fields}, Grundlehren Math. Wiss.
\textbf{323}, Springer (2000).


\bibitem{RW} J.\ Ritter and A.\ Weiss, Towards equivariant
Iwasawa theory II, \textit{Indag. Math. (N.S.)}, 15(4) (2004),
549-572.

\bibitem{RZ} L.\ Ribes and P.\ Zalesskii, \textit{Profinite Groups}, Ergeb. Math. Grenzgeb.
\textbf{3}, Springer-Verlag, Berlin, (2000).


\bibitem{Wei} C.\ A.\ Weibel, \textit{An Introduction to
Homological Algebra}, Reprinted, Cambridge Stud. Adv. Math.
\textbf{38}, Cambridge Univ. Press, Cambridge, UK, (1997).


\bibitem{Wil} J.\ Wilson, \textit{Profinite Groups}, London Mathematical
Society Monographs New Series, vol. \textbf{19}, Oxford University
Press, (1998).

\end{thebibliography}
\end{document}